\RequirePackage{fix-cm}
\let\accentvec\vec

\documentclass[smallcondensed]{svjour3}     
\smartqed  
\usepackage{graphicx}
 \journalname{}

\usepackage{amssymb}

\let\vec\accentvec

\usepackage{amsmath}

\usepackage{comment}
\usepackage{enumerate}

\usepackage{hyperref}

\newcommand{\rd}{{\mathbb R}^d}
\newcommand{\zd}{{\mathbb Z}^d}
\newcommand{\td}{{\mathbb T}^d}
\newcommand{\n}{{\mathbb N}}
\newcommand{\z} {{\mathbb Z}}
\newcommand{\rn}{{\mathbb R}}
\newcommand{\cn} {{\mathbb C}}
\newcommand{\w} {\widetilde}
\newcommand{\h} {\widehat}
\newcommand{\nul} {{\bf 0}}

\begin{document}

\title{A directional uncertainty principle for periodic functions\thanks{The work is supported by Volkswagen Foundation; the first and the second authors are also
supported by grants from RFBR \# 15-01-05796-a,  St. Petersburg State University \#~9.38.198.2015.}
\thanks{This is a pre-print of an article submitted to Multidimensional Systems and Signal Processing. The final authenticated version is available online at:
\\ \url{https://doi.org/10.1007/s11045-018-0613-1}
\\ \url{https://rdcu.be/49OH}}
}


\author{A. Krivoshein         \and
        E. Lebedeva \and 
        J. Prestin 
}


\institute{A. Krivoshein  \at
              St. Petersburg State University
              \email{krivosheinav@gmail.com}           
  \and
           E. Lebedeva \at
             Peter the Great St.Petersburg Polytechnic University\at
              St. Petersburg State University
              \email{ealebedeva2004@gmail.com }  
          \and     
           J. Prestin  \at
              University of L\"ubeck
              \email{prestin@math.uni-luebeck.de}  
}

\date{Received: date / Accepted: date}

\maketitle

\begin{abstract}
In this paper we introduce a notion of a directional uncertainty product for multivariate periodic functions and multivariate discrete signals. It measures a 
localization of a signal along a particular direction. We study properties of the uncertainty product and give an example of well localized multivariate periodic Parseval wavelet frames.
\end{abstract}

\section{Introduction}
\label{intro}
A notion of uncertainty product is a sufficiently well-studied object in harmonic analysis. 
Initially, it was introduced for functions 
on the real line
to measure a simultaneous localization of a function and its Fourier transform~\cite{Schr}. The essence of this measurement is concentrated in the fundamental Heisenberg uncertainty principle, which says that for any appropriate function the uncertainty product cannot be smaller than a positive absolute constant.     
Later, numerous versions of this framework were developed for different algebraic and topological
structures such as abstract  locally compact groups, high-dimensional  spheres, etc. (see, e.g.,~\cite{KL15},~\cite{NW},~\cite{PS}). For more detailed information concerning this topic, we refer the interested reader to surveys~\cite{FS} and~\cite{Ricaud2014} and the
references therein.

In this paper we focus on the case of multivariate periodic functions and multivariate discrete signals. For periodic functions of one variable a 
 notion of uncertainty product
 was introduced in 1985 by Breitenberger in~\cite{B}. The corresponding uncertainty principle is also valid in this setup.
One possible extension of this notion to the case of multivariate periodic functions was suggested by Goh and Goodman in~\cite{GohGoodman2004UC} (see formula (\ref{eq:UCGGdef})). However, this approach does not take into account the main difference between periodic functions of one variable and many variables, namely the localization of a function along particular directions. The main contribution of this paper is a new approach  that allows to include the directionality into the definition of the uncertainty product (see formula (\ref{eq:UCLdef})). We compare these two approaches and show that they are not equivalent (see Lemma \ref{lem:Compare}). 
At the same time, both definitions fit into a more general operator approach (see formula (\ref{eq:UCmain})). This approach was established by Folland   in~\cite{Folland1989Book} and was extended to two normal or symmetric operators by Selig in~\cite{Selig2002UC}  and Goh, Micchelli in~\cite{GohMicch2002UC}. For several operators this approach was generalized by Goh and Goodman in~\cite{GohGoodman2004UC}.

From the other point of view, this directional uncertainty product is applicable for multidimensional discrete signals due to the duality: periodic signal - discrete spectrum (Fourier series) and discrete signal - periodic spectrum (the Discrete-Time Fourier Transform, DTFT). In this sense, our definition is an alternative to the one given in~\cite{Meerkoetter2017} and allows to take into account the directionality of signals.


The paper is organized as follows. Section 2 is devoted to basic definitions.
In Section 3 we study  the properties of the directional uncertainty product  for periodic functions  and compare this product with one defined by Goh and Goodman.
Lemma \ref{lem:UCTdOptimalSeq} gives a sequence of trigonometric polynomials such that the sequence of their
directional uncertainty products tends to the optimal value.
Lemma~\ref{lem:Compare} illustrates a difference between these two uncertainty products. 
In   Subsection 3.1. we study the behavior of both uncertainty products for the Dirichlet and Fej\'er kernels. Lemmas~\ref{lem:DirichletUCTD} and~\ref{lem:FejerUCTD} concern the directional case. In Lemma~\ref{lem:DFUCTD} we address the same question to the Goh and Goodman case. In Subsection 3.2. we minimize  the directional angular  variance  for
trigonometric polynomials. Theorem~\ref{theo:MinAngL} describes the case of the directional uncertainty product, and Theorem~\ref{theo:MinAngGG} corresponds to the case of the uncertainty product defined by Goh and Goodman. 
In Section 4 we give an example of a multivariate periodic Parseval wavelet frame with a small directional uncertainty product (see Theorem \ref{main_loc}).

\section{Basic notations and definitions}
\label{sec2}
We use the standard multi-index notation.
    Let $d\in\n,$ $\rd$ be the $d$-dimensional Euclidean space, $\{e_j, 1\le j\le d\}$ be the standard basis in $\rd$, 
     $\zd$ be the integer lattice  in $\rd$, 
    $\td=\rd\slash\zd$ be the $d$-dimensional torus.
    Let  $x = (x_1,\dots, x_d)^\mathrm{T}$ and 
    $y =(y_1,\dots, y_d)^\mathrm{T}$ be column vectors in $\rd$.
    Then $\langle x, y\rangle:=x_1y_1+\dots+x_dy_d$,
    $\|x\| := \sqrt {\langle x, x\rangle}$, $\|x\|_1 = \sum_{j=1}^d |x_j|$, $\|x\|_{\infty} = \max_j |x_j|$.
    We say that
     $x \geq y,$ if $x_j \geq y_j$ for all $j= 1,\dots, d,$
    and we say that 
 	$x > y,$ if $x \geq y $ and $x\neq y.$ 
    Further, $\zd_+:=\{\alpha\in\zd:~\alpha\geq~{\bf0}\},$
    where  ${\bf0}=(0,\dots, 0)$ denotes the origin in $\rd$.
    For $\alpha=(\alpha_1,\dots,\alpha_d)^\mathrm{T}\in \zd_+$,
    denote $|\alpha|:=\alpha_1+ \dots +\alpha_d.$ 
    For $x\in {\mathbb R}$, $x_+ := \begin{cases}
    0, & x\le 0,\\
    x, & x>0.
\end{cases}$
    
	For a sufficiently smooth  function $f$ defined on $\Omega\subset\rd$ and a multi-index $\alpha\in\zd_+$, 
	$D^{\alpha} f$ denotes the derivative of $f$ of order $\alpha$ and
	$D^{\alpha} f =\frac{\partial^{|\alpha|} f}
    {\partial x^{\alpha}}=\frac{\partial^{|\alpha|} f
    }{\partial^{\alpha_1}x_1\dots
    \partial^{\alpha_d}x_d}.$ For $\alpha = e_j$, we also use $D^{e_j} f = f'_j.$
	The directional derivative of a sufficiently smooth  function $f$ defined on $\Omega$ along a vector $L=(L_1,...,L_d)\in\rd$ is denoted by
	$\frac{\partial f}{\partial L} =\sum_{j=1}^d L_j \frac{\partial f}{\partial x_j}.
	$

For a function $f\in L_2(\td)$ its norm is denoted by 
    $\|f\|^2_{\td} = \int_{\td} |f(x)|^2 \mathrm{d}x$.
    The Fourier series coefficients of a function $f\in L_2(\td)$ are given by $c_k = c_k(f) = \h f(k) = \int_{\td} f(x) \mathrm{e}^{-2\pi \mathrm{i} \langle k,x\rangle} \mathrm{d}x$, $k\in\zd.$
    The Sobolev space $H^1(\td)$ consists of functions in $L_2(\td)$ such that all its derivatives of the first order are also in $L_2(\td)$, which can be written as
    $$
    H^1(\td)= \left\{ f\in L_2(\td): \sum_{k\in\zd} \|k\|^2 |c_k(f)|^2 < \infty \right\}.
    $$

Let ${\cal H}$ be a Hilbert space with inner product $\langle \cdot, \cdot \rangle$ and with norm $\|\cdot\| := \langle \cdot, \cdot \rangle^{1/2}$. 
Let ${\cal A}$, ${\cal B}$ be two linear operators with domains ${\cal D}({\cal A})$, ${\cal D}({\cal B})\subseteq {\cal H}$ and ranges in ${\cal H}$.
The variance of non-zero $f\in {\cal D}({\cal A})$ with respect to the operator ${\cal A}$ is defined to be
 $$\Delta_{}({\cal A},f) = \|{\cal A} f\|^2 - \frac{|\langle {\cal A} f, f\rangle|^2}{\|f\|^2} = \left\|\left({\cal A} - \frac{\langle {\cal A} f, f\rangle}{\|f\|^2}\right) f \right\|^2.$$ 
The commutator of ${\cal A}$ and ${\cal B}$ is defined  by 
$
[{\cal A},{\cal B}] := {\cal A}{\cal B}-{\cal B}{\cal A}
$
with domain $ {\cal D}({\cal A}{\cal B}) \bigcap {\cal D}({\cal B}{\cal A})$.

\begin{theorem}~\cite[Theorem 4.1]{GohGoodman2004UC}
    Let ${\cal A}_1,\dots{\cal A}_n,$ ${\cal B}_1,\dots{\cal B}_n$ be symmetric or normal operators acting from a Hilbert space ${\cal H}$ into itself. Then for any non-zero $f$ in ${\cal D}({\cal A}_j{\cal B}_j) \bigcap {\cal D}({\cal B}_j{\cal A}_j)$, $j=1,\dots,n$,
    \begin{equation}
    \frac{1}{4} \left( \sum_{j=1}^n |\langle [{\cal A}_j,{\cal B}_j]  f, f\rangle| \right)^2 \le \left(\sum_{j=1}^n \Delta_{}({\cal A}_j,f)\right)\left(\sum_{j=1}^n \Delta_{}({\cal B}_j,f)\right).
    \label{eq:UCmain}
    \end{equation}
\end{theorem}

If the commutator $\langle [{\cal A}_j,{\cal B}_j]  f, f\rangle$  is non-zero for all $j=1,\dots,n$, then the uncertainty product for $f$ is defined as
$$
\mathrm{UP}(f):= \left(\sum_{j=1}^n \Delta_{}({\cal A}_j,f)\right)\left(\sum_{j=1}^n \Delta_{}({\cal B}_j,f)\right) \left( \sum_{j=1}^n |\langle [{\cal A}_j,{\cal B}_j]  f, f\rangle| \right)^{-2} .
$$
In this terms, the uncertainty principle says that the uncertainty product $\mathrm{UP}(f)$ cannot be smaller than $\frac{1}{4},$ for any appropriate function $f$. 

The well-known Heisenberg uncertainty product for functions defined on the real line fits in this operator approach, if $n=1$, ${\cal H} = L_2({\mathbb R})$ and the two operators are as follows ${\cal A} f(x)= 2 \pi x f(x)$, ${\cal B} f(x)= \frac{\mathrm{i} }{2\pi}\frac{\mathrm{d}f}{\mathrm{d}x}(x)$. 

The Breitenberger uncertainty product is defined for  periodic functions. In this case, $n=1$, ${\cal H} = L_2({\mathbb T})$  and ${\cal A}^{{\mathbb T}} f(x)= \mathrm{e}^{2\pi \mathrm{i} x} f(x)$, ${\cal B}^{{\mathbb T}} f(x)= \frac{ \mathrm{i} }{2\pi}\frac{\mathrm{d}f}{\mathrm{d}x}(x)$. The commutator is $[{\cal A}^{{\mathbb T}},{\cal B}^{{\mathbb T}}]={\cal A}^{{\mathbb T}}$. It is more convenient for the Breitenberger uncertainty product, to use the notions of the angular and frequency variance. Since $\|{\cal A}^{{\mathbb T}} f\|^2_{{\mathbb T}} = \|f\|^2_{{\mathbb T}}$,
$$\text{var}^A(f) = \frac{\|f\|^2_{{\mathbb T}} \Delta({\cal A}^{{\mathbb T}},f)}{|\langle [{\cal A}^{{\mathbb T}},{\cal B}^{{\mathbb T}}] f,f\rangle|^2} =  \left(\frac{\|f\|^2_{{\mathbb T}}}{|\langle{\cal A}^{{\mathbb T}} f,f\rangle|}\right)^2 - 1 , $$
   $$ \text{var}^F(f) = \frac{\Delta({\cal B}^{{\mathbb T}},f)}{\|f\|^2_{{\mathbb T}}} = \frac{\|{\cal B}^{{\mathbb T}} f\|^2_{\rd}}{\|f\|^2_{{\mathbb T}}} - \frac{|\langle{\cal B}^{{\mathbb T}} f, f\rangle|^2}{\|f\|^4_{{\mathbb T}}},$$
    $$
  \mathrm{UP}^{{\mathbb T}}(f) := \text{var}^{A}(f) \text{var}^F(f).
  $$
It is known that the lower bound for $\mathrm{UP}^{{\mathbb T}}$ does not attain on any function. But there exist sequences of functions such that $\mathrm{UP}^{{\mathbb T}}$ tends to the optimal value $\frac{1}{4}$ (see, e.g.,~\cite{PrestinQ99UC}).


For the space $L_2(\td)$ of multivariate periodic functions, Goh and Goodman in~\cite{GohGoodman2004UC} suggest to take the operators as follows ${\cal A}_j f(x) = \mathrm{e}^{2\pi \mathrm{i} x_j} f(x)$, ${\cal B}_j f(x)= \frac{ \mathrm{i} }{2\pi} \frac{\partial f}{\partial x_j}(x)$, $j=1,\dots,d.$ Note that the domains of the operators are $\bigcap_{j=1}^d {\cal D}({\cal A}_j) = L_2(\td),$ $\bigcap_{j=1}^d {\cal D}({\cal B}_j) = H^1(\td)$. Operators ${\cal A}_j$ are normal, ${\cal B}_j$ are self-adjoint.
The commutators for $f\in H^1(\td)$ are $[{\cal A}_j,{\cal B}_j]f=  {\cal A}_j f.$ 
The uncertainty principle for these operators is stated as follows.
\begin{theorem}
    For a function $f\in H^1(\td)$, such that 
    $\langle {\cal A}_j f, f\rangle\neq 0$, $j=1,\dots,d,$ the uncertainty product $\mathrm{UP}^{\td}_{GG}(f)$ is well-defined and
{\small
\begin{equation}
    \mathrm{UP}_{GG}^{\td}(f)= \frac{\sum\limits_{j=1}^d \left(\|f\|_{\td}^4 -\left|\sum\limits_{k\in\zd} c_{k-e_j}\overline{c_k}\right|^2\right)}{\left(\sum\limits_{j=1}^d\left|\sum\limits_{k\in\zd} c_{k-e_j}\overline{c_k}\right|\right)^2}
\sum\limits_{j=1}^d\left( \frac{\sum\limits_{k\in\zd} k_j^2 |c_k|^2}{\|f\|_{\td}^2} - \left(\frac{\sum\limits_{k\in\zd} k_j|c_k|^2}{\|f\|_{\td}^2} \right)^2 \right) \ge \frac 14,
\label{eq:UCGGdef}
\end{equation}
}
where $k=(k_1,\dots,k_d)$, $c_k=c_k(f)$ are the Fourier coefficients of $f.$ 
\end{theorem}

\noindent Defining the variances for $f\in H^1(\td)$ as $$\text{var}_{GG}^A(f) = \frac{\|f\|^2_{{\mathbb T}^d} \sum\limits_{j=1}^d \Delta({\cal A}_j,f)}{\left( \sum_{j=1}^n |\langle [{\cal A}_j,{\cal B}_j]  f, f\rangle| \right)^2}, \quad \text{var}_{GG}^F(f) = \sum\limits_{j=1}^d \Delta({\cal B}_j,f)/\|f\|^2_{{\mathbb T}^d},$$ it can be shown, that the variances attain the value $\infty$ if and only if $\langle {\cal A}_j f, f\rangle= 0$, for all $j=1,\dots,d.$ 
In these cases, we can also assign to $\mathrm{UP}_{GG}^{\td}(f)$ the value $\infty$, except the following case $\text{var}_{GG}^F(f) = 0$ and $\text{var}_{GG}^A(f) = \infty$. This case happens  if and only if $f$ is a monomial, since $\text{var}_{GG}^F(f) = 0$ if and only if $f$ is a monomial.
Indeed, since $\text{var}_{GG}^F(f)=\text{var}_{GG}^F(a f)$ for any appropriate $f\in H^1(\td)$ and $a\in\rd,$ $a\neq 0,$ we can assume that $\|f\|_{\td} = 1$. Therefore, 
$$
    \text{var}_{GG}^F(f) = 0
    \quad \text{ if and only if }
     \left| \left\langle {\cal B}_j f, f\right\rangle \right|^2 =  \left\|  {\cal B}_j f \right\|_{\td}^2 \left\| f \right\|_{\td}^2 \quad \forall j=1,\dots,d.
$$
Due to the Cauchy-Bunyakovsky-Schwarz inequality the equality is possible only if
$ f = \alpha_j \frac{\partial f}{\partial x_j}$, where $\alpha_j \in \cn$ for all $j=1,\dots,d.$
Thus, $f$ should be a monomial.
However, in this case, i.e., $\text{var}_{GG}^F(f) = 0$ and $\text{var}_{GG}^A(f) = \infty$,
inequality~(\ref{eq:UCmain}) takes the form 
$1/4 \cdot 0 \le C \cdot 0.$ It is trivially true. Thus,  inequality~(\ref{eq:UCmain}) is valid for all non-zero functions $f\in H^1(\td).$ 

In fact, the above approach for the definition of the uncertainty product does not deal with a new phenomenon, that appears in the multidimensional case, namely, the localization of a function along particular directions. We suggest an approach that allows to include the directionality into the definition.

The directional uncertainty product for $\td$ along a direction $L\in\zd$ ($L\neq \nul$) is defined using the operators  
  $${\cal A}_Lf(x) = \mathrm{e}^{2\pi \mathrm{i} \langle L, x\rangle} f(x), \quad {\cal B}_L f(x) = \frac{ \mathrm{i} }{2\pi} \frac{\partial f}{\partial L}(x).$$
with domains 
${\cal D}({\cal A}_L) = L_2(\td),$ ${\cal D}({\cal B}_L) = H^1(\td)$. Note that
 ${\cal A}_L$ is normal, ${\cal B}_L$ is self-adjoint. 
The commutator for $f\in {\cal D}({\cal A}_L) \cap {\cal D}({\cal B}_L)$ is $      [{\cal A}_L,{\cal B}_L]f= \|L\|^2 {\cal A}_L f.$
Thus, the directional uncertainty product for a function $f\in {\cal D}({\cal A}_L) \cap {\cal D}({\cal B}_L)$ such that $ {\cal A}_L f \neq 0$ is defined as 
    {\small  $$
  \mathrm{UP}^{\td}_L(f) = \frac{1}{\|L\|^4_2}\left( \frac{\|f\|^4_{\td}}{|\langle{\cal A}_Lf,f\rangle|^2} - 1 \right)\left( \frac{\|{\cal B}_L f\|^2_{\td}}{\|f\|^2_{\td}} - \frac{|\langle{\cal B}_L f, f\rangle|^2}{\|f\|^4_{\td}} \right):= \frac{1}{\|L\|^4}\text{var}_L^A(f) \text{var}_L^F(f),
  $$ }
  where $\text{var}_L^A(f)$ is the angular directional variance and $\text{var}_L^F(f)$ is the frequency directional variance.

\begin{theorem}
    For $L\in\zd$ and a function $f\in H^1(\td)$, such that 
    $\langle {\cal A}_L f, f\rangle\neq 0$, the uncertainty product $\mathrm{UP}^{\td}_{L}(f)$ is well-defined and
    {\small
    \begin{equation}
\mathrm{UP}^{\td}_L(f) = \frac{1}{\|L\|^4} \left( \frac{\left(\sum\limits_{k\in\zd} |c_k|^2\right)^2}{\left|\sum\limits_{k\in\zd} c_{k-L}\overline{c_k}\right|^2} - 1 \right)
\left( \frac{\sum\limits_{k\in\zd} 
\langle L, k\rangle^2 |c_k|^2}{\sum\limits_{k\in\zd} |c_k|^2} - \left(\frac{\sum\limits_{k\in\zd} \langle L, k\rangle |c_k|^2}{\sum\limits_{k\in\zd} |c_k|^2} \right)^2 \right) \ge \frac 14,
\label{eq:UCLdef}
\end{equation}
}
where $c_k=c_k(f)$ are the Fourier coefficients of $f.$ 
\end{theorem}
The statement  easily follows from the operator approach and 
   $${\cal A}_L f(x) = \sum\limits_{k\in \zd} c_{k-L} \mathrm{e}^{2 \pi \mathrm{i} \langle k, x\rangle}, \quad
   {\cal B}_L f(x) = - \sum\limits_{k\in \zd} \langle L, k\rangle c_{k} \mathrm{e}^{2 \pi \mathrm{i} \langle k, x\rangle}.$$
    
It can be shown, that the directional variances attain the value $\infty$ if and only if $\langle {\cal A}_L f, f\rangle= 0$. 
In this case, we can also assign to $\mathrm{UP}_{L}^{\td}(f)$ the value $\infty$, except the following case $\text{var}_{L}^F(f) = 0$ and $\text{var}_{L}^A(f) = \infty$. This case happens  if and only if $f$ is a monomial (with the arguments as above).
Thus analogously,  inequality~(\ref{eq:UCmain}) for operators ${\cal A}_L$ and ${\cal B}_L$ is valid for all non-zero functions $f\in H^1(\td).$



\section{The properties of the directional uncertainty product for the periodic case}

First of all, we note that the standard manipulations of functions like shifts, modulations and multiplying by numbers do not change the uncertainty product $\mathrm{UP}_L^{\td}.$ 
  
\begin{lemma}
    Let $f\in H^1(\td)$.
    Suppose $g(x) = a \ \mathrm{e}^{2\pi \mathrm{i}\langle K, x\rangle} f(x-x_0),$ where $K\in\zd$, $a\in\mathbb{R}$, $a\neq 0$, $x_0\in\rd,$  then 
    $\mathrm{UP}^{\td}_L(g) = \mathrm{UP}^{\td}_L(f).$
\end{lemma}

The proof can be done by straightforward computations.

As for the Breitenberger uncertainty product and for the uncertainty product defined by Goh and Goodman,  the  optimal function for the directional uncertainty product does not exist.
Indeed, let $a(f)= \frac{\langle{\cal A}_L f,f\rangle}{\|f\|_2^2}$ and
$b(f)= \frac{\langle{\cal B}_L f,f\rangle}{\|f\|_2^2}$. Since ${\cal B}_L$ is self-adjoint, $b(f)$ is real.
Due to Theorem 3.1 in~\cite{Selig2002UC} the equality for the uncertainty principle is attained if and only if there exist $\lambda\in{\mathbb C}$ such that
$$
({\cal B}_L - b(f))f = \lambda ({\cal A}_L - a(f))f = -\overline{\lambda} ({\cal A}^*_L - \overline{a(f)}) f.
$$
The second identity yields
\begin{multline*}
f(x) \left( \lambda \mathrm{e}^{2\pi \mathrm{i} \langle L, x\rangle} + \overline{\lambda} \mathrm{e}^{-2\pi \mathrm{i} \langle L, x\rangle} - a(f) \lambda - \overline{\lambda} \overline{a(f)}\right) \\ = 2 f(x) (\mathcal{R}e(\lambda \mathrm{e}^{2\pi \mathrm{i} \langle L, x\rangle}) - \mathcal{R}e(a(f) \lambda)) \equiv 0.
\end{multline*}
This condition can be satisfied only if $f = 0$ or $\lambda=0.$ For the second case, we get $({\cal B}_L - b(f))f 
=0$ or $\frac{ \mathrm{i} }{2\pi} \frac{\partial f}{\partial L}(x) = b(f) f(x),$
which is only possible when $f$ is a monomial, i.e. $f(x) = C \mathrm{e}^{2\pi \mathrm{i} \langle k, \xi\rangle}$.
Recall that for monomials the directional uncertainty product is not defined.

The next lemma gives a sequence of trigonometric polynomials such that the sequence of their directional uncertainty products tends to the optimal value.

\begin{lemma}
    Suppose $p_n(x)=(1+\cos 2\pi \langle L, x\rangle)^n$ for $n\in\n.$ Then
    $$
    \mathrm{UP}^{\td}_L(p_n) = \frac{1}{4} + O\left(\frac{1}{n}\right), \quad \text{ as } n\to\infty.
    $$
    \label{lem:UCTdOptimalSeq}
\end{lemma}

{\bf Proof.}
Denote 
    $I_{2n} := \int_{\td} (1+\cos 2\pi \langle L, x\rangle)^{2n} \mathrm{d}x = \|p_n\|^2_{\td}.$
    Since $p_n$ is even 
$\langle {\cal A}_L p_n, p_n\rangle   = I_{2n+1} - I_{2n}.$
Further,
\begin{equation*}
    \begin{split}
    \|{\cal B}_L p_n\|^2_{\td}  & = 
    n^2 \|L\|^4 \int_{\td} (1+\cos 2\pi \langle L, x\rangle)^{2n-2} \sin^2 2\pi \langle L, x\rangle \ \mathrm{d}x \\ & = n^2 \|L\|^4 (2I_{2n-1} - I_{2n}),
    \end{split}
\end{equation*}
    since $\sin^2 2\pi \langle L, x\rangle = 2 (1+\cos2\pi \langle L, x\rangle)- (1+\cos 2\pi \langle L, x\rangle)^2.$
    Again, since $p_n$ is even and $\sin 2\pi \langle L, x\rangle$ is odd, we get
    $\langle {\cal B}_L p_n, p_n\rangle  =0.$
    So, finally, 
    \begin{equation}
        \begin{split}          \mathrm{UP}^{\td}_L(p_n)
    & = n^2  \frac{I_{2n+1}}{I_{2n}}
    \frac{(2I_{2n} - I_{2n+1})(2I_{2n-1} - I_{2n})}{(I_{2n+1} - I_{2n})^2}.
        \end{split}
        \label{eq:UCpn}
    \end{equation}
    It remains to compute $I_n.$ Since
    $(1+\cos 2\pi \langle L, x\rangle)^{n}=2^n \cos^{2n} \frac{ 2\pi \langle L, x\rangle}{2},
    $
$$
    I_n 
    = 
    \frac{2^n}{2^d} \int\limits_{[0,2)^d} \cos^{2n} \frac{ 2\pi \langle L, x\rangle}{2} \mathrm{d}x = 
    2^n \int\limits_{[0,1)^d} \cos^{2n} (2\pi \langle L, x\rangle) \mathrm{d}x.
$$
Using Euler's formula we get
\begin{equation*}
    \begin{split}
    \cos^{n} (2\pi \langle L, x\rangle) & = \frac{1}{2^n} \sum\limits_{j=0}^n \binom{n}{j} \mathrm{e}^{2\pi \mathrm{i} j \langle L, x\rangle}\mathrm{e}^{-2\pi \mathrm{i} (n-j) \langle L, x\rangle}  \\ & = 
    \frac{1}{2^n} \sum\limits_{j=0}^n \binom{n}{j} \mathrm{e}^{2\pi \mathrm{i} (2j-n) \langle L, x\rangle}.
    \end{split}
\end{equation*}
Due to the Parseval equality for the function $\cos^{n} (2\pi \langle L, x\rangle)$ we obtain
    $$
   \int_{\td} \cos^{2n} (2\pi \langle L, x\rangle) \mathrm{d}x = 
    \frac{1}{2^{2n}} \sum\limits_{j=0}^n \binom{n}{j}^2 =   \frac{1}{2^{2n}}\binom{2n}{n} = \frac{1}{2^{2n}}\frac{(2n)!}{(n!)^2}.
    $$
Therefore,    $ I_n = \frac{(2n-1)!!}{n!}.$
 Here $(2n-1)!!$ is the double factorial of $2n-1$.
Substituting this    in~(\ref{eq:UCpn}), we obtain
    $
    \mathrm{UP}^{\td}_L(p_n) = 
    \frac{1}{4}  + \frac{1}{8n-2}.      $ $\Diamond$

Let us compare the uncertainty product defined by Goh and Goodman and the directional uncertainty product.
They are not equivalent. The next lemma gives a pair of examples where the uncertainty products behave differently.


\begin{lemma}
    Let $L\in\zd$. 
    \begin{enumerate}[(A)]
        \item Suppose $\w p_n(x)=(1+\cos 2\pi \langle L, x\rangle)^n+2 \cos 2 \pi x_1,$  where
 $L$ is not collinear to $e_1$. Then
$$ \mathrm{UP}^{\td}_L(\w p_n) \to \frac{1}{4},\qquad \frac{\mathrm{UP}^{\td}_{GG}(\w p_n)}{n \ 4^n} \to  \frac{d \|L\|^2}{32}  \quad n\to\infty. $$
\item Suppose $\w t_n(x)=(1+\cos 2 \pi x_1)^n+2 \cos 2 \pi \langle L, x\rangle$,  where
 $L$ is not collinear to all $e_j$ and $|L_j|>1$ for all $j=1,...,d$. Then
$$ \frac{\mathrm{UP}^{\td}_{L}(\w t_n)}{n \ 4^n} \to \frac{L_1^2}{32 \|L\|^4},\qquad \frac{\mathrm{UP}^{\td}_{GG}(\w t_n)}{n} \to \frac{d-1}{4},  \quad n\to\infty. $$
    \end{enumerate}
    \label{lem:Compare}
\end{lemma}

{\bf Proof.} Let us prove item (A).
For convenience, we will use the
notation $p_n(x)=(1+\cos 2\pi \langle L, x\rangle)^n$ and some facts used in the proof of Lemma~\ref{lem:UCTdOptimalSeq}. 
Then
$$
\|\w p_n\|^2_{\td} = \|p_n\|^2_{\td}  + 2 =  \frac{(4n-1)!!}{(2n)!} + 2, \qquad
\langle {\cal A}_L \w p_n, \w p_n \rangle=\langle {\cal A}_L  p_n,  p_n \rangle = I_{2n+1} - I_{2n},
$$
    $$
    {\cal B}_L \w p_n(x) =  i  \|L\|^2 n (1+\cos (2\pi \langle L, x\rangle))^{n-1} \sin (2\pi \langle L, x\rangle) + 2 i L_1 \sin 2 \pi x_1,
    $$
    $$
    \|{\cal B}_L \w p_n\|^2_{\td} =
    n^2 \|L\|^4 (2I_{2n-1} - I_{2n}) + 2 L_1^2.
    $$
Since $\w p_n$ is even  and  ${\cal B}_L  \w p_n$ is odd we get
$\langle{\cal B}_L  \w p_n, \w p_n \rangle  =0.$
Therefore, 
{\footnotesize
 $$\mathrm{UP}^{\td}_L(\w p_n) =  \frac{1}{ \|L\|^4 }  \left( \frac{ \left(\frac{(4n-1)!!}{(2n)!} + 2\right)^2}{\left(2n\frac{ (4n-1)!!}{(2n+1)!}\right)^2} -1 \right) 
 \left( \frac{ n^2 \|L\|^4 \frac{(4n-3)!!}{(2n)!}  + 2 L_1^2 }{\frac{(4n-1)!!}{(2n)!} + 2} \right) $$  
   $$
=\frac{n^2}{(2n+1)(4n-1)}  \frac{ \left(1 + 2\frac{(2n)!(2n+1)}{(4n-1)!!} \right)\left(2 + 2\frac{(2n)!}{(4n-1)!!} - \frac{1}{2n+1} \right)}{\left(\frac{2n}{2n+1}\right)^2} 
 \left( \frac{1  + 2 \frac{L_1^2}{\|L\|^4} \frac{(2n)! (4n-1)}{n^2 (4n-1)!!} }{1  + 2 \frac{(2n)!}{(4n-1)!!}} \right).
 $$}
By the Stirling formula $n! = \sqrt{2 \pi n} \left( \frac{n}{\mathrm{e}}\right)^n (1 + O(1/n))$, it follows that 
$\frac{(2n)!}{(4n-1)!!} = \frac{\sqrt{2 \pi n} (1 + O(\frac{1}{n})) }{2^{2n} } \to 0,$  $n\to\infty.$
Therefore,
$ \mathrm{UP}^{\td}_L(\w p_n) \to \frac{1}{4}, \quad n\to\infty.$

Now, we compute $\mathrm{UP}^{\td}_{GG}(\w p_n)$. 
Let $\w c_k = \w c_k(\w p_n)$ be the Fourier coefficients of $\w p_n.$ Then 
$$
\w c_0 = \int_{\td} \w p_n(x) \mathrm{d}x = \int_{\td} p_n(x) \mathrm{d}x = I_n = \frac{(2n-1)!!}{n!},
$$
$$
\langle {\cal A}_j \w p_n, \w p_n \rangle = \sum\limits_{k\in\zd} \w c_{k-e_j} \w c_k = 
\delta_{j,1} (\w c_{e_1} \w c_{\nul} + \w c_\nul \w c_{e_1}) = 
2 \delta_{j,1} \frac{(2n-1)!!}{n!}, 
$$
for $j=1,\dots,d.$ Further,
    $$
    B_j \w p_n(x) 
    = - \mathrm{i} L_j n (1+\cos (2\pi \langle L, x\rangle))^{n-1} \sin (2\pi \langle L, x\rangle) - 2 \mathrm{i} \delta_{j,1} \sin 2 \pi x_1.
    $$
Therefore,
    $  \|B_j \w p_n\|^2_{\td} = 
     n^2 L_j^2  (2I_{2n-1} - I_{2n}) +2 \delta_{j,1}.$
Since   $\w p_n$ is even  and $B_j \w p_n$ is odd, we get $
    \langle B_j \w p_n, \w p_n \rangle  =0.$
Hence, combining all results in the definition of $\mathrm{UP}^{\td}_{GG}(\w p_n)$~(\ref{eq:UCGGdef}) and after some simplifications, we obtain
{\footnotesize
 $$\mathrm{UP}^{\td}_{GG}(\w p_n) =  
 \frac{n^2\|L\|^2}{4(4n-1)}  \left( d  \left(\frac{(4n-1)!!}{(2n)!} \frac{n!}{(2n-1)!!} + 2 \frac{n!}{(2n-1)!!} \right)^2 - 4  \right) 
\frac{1 +  \frac{2 (2n)!}{n^2 \|L\|^2 (4n-1)!!} }{1+ 2\frac{(2n)!}{(4n-1)!!}}.$$ }
By the Stirling formula $\frac{(2n)!}{(4n-1)!!} =  \frac{\sqrt{2 \pi n} (1 + O(\frac{1}{n})) }{2^{2n} } \to 0$  as $n\to\infty$ and
$\frac{n!}{(2n-1)!!} = \frac{\sqrt{ \pi n} (1 + O(\frac{1}{n})) }{2^{n} } \to 0$ as $n\to\infty$. Thus, 
$\frac{(4n-1)!!}{(2n)!} \frac{n!}{(2n-1)!!} = \frac{2^n}{\sqrt{2}} (1 + O(\frac{1}{n}))$ as $n\to\infty$.
Finally, it follows that
$\frac{\mathrm{UP}^{\td}_{GG}(\w p_n)}{n 4^n} \to  \frac{d \|L\|^2}{32}$ as $n\to\infty$.

Item (B) can be proved analogously. By similar arguments it can be shown that
{\footnotesize
 $$\mathrm{UP}^{\td}_L(\w t_n) = 
   \frac{1}{ \|L\|^4 }  \left(  \left(\frac{\frac{(4n-1)!!}{(2n)!} }{\frac{ (2n-1)!!}{n!}} \right)^2 \frac{ \left(1 + 2\frac{(2n)!}{(4n-1)!!}\right)^2}{4} -1 \right) 
\frac{L_1^2/2  +  2 \|L\|^4 \frac{(2n-1)!} {n (4n-3)!!} }{ 1+ 2\frac{(2n)!}{(4n-1)!!}} \frac{2n^2}{4n-1}
 $$}
and
{\footnotesize
$$
 \mathrm{UP}^{\td}_{GG}(\w t_n) = \left( d \left(\frac{2n+1}{2n} +2 \frac{(2n)!}{(4n-1)!!} \frac{2n+1}{2n}\right)^2 - 1 \right)
\frac{\frac{n}{2 } + 2  \|L\|^2 \frac{(2n-1)!}{(4n-3)!!}}{1 + 2 \frac{(2n)!}{(4n-1)!!}} \frac{2n}{4n-1}.
$$}
The Stirling formula yields Item (B).$\Diamond$

\subsection{The uncertainty products for the Dirichlet and Fej\'er kernels}
  
As it was noted in~\cite{PrestinQ99UC}, the sequence of the Breitenberger uncertainty products of the Dirichlet kernels $D_n(x)=\sum_{k=-n}^n \mathrm{e}^{2\pi \mathrm{i} k x}$ tends to infinity as $n\to \infty.$ In~\cite{PrestinQ95TimeFreq} it was noted that the sequence of Breitenberger uncertainty products of the Fej\'er kernels $F_n(x)=\sum\limits_{k=-n}^n (1-|k|/n)\mathrm{e}^{2\pi \mathrm{i} k x}$ tends to $\frac{3}{10}$ as $n\to \infty.$ In the multivariate case the analogous difference between these kernels also holds for the directional uncertainty product and the one defined by Goh and Goodman. Different methods of summation can be used for the Dirichlet kernel. Let us consider a rectangular one. 
\begin{lemma}
    Let $
D_N(x)= \sum\limits_{-N \le k \le N} \mathrm{e}^{2\pi \mathrm{i} \langle k, x\rangle},
$
where $N\in \zd,$ $N>\nul,$ $L\in\zd.$ Then 
$$ \mathrm{UP}^{\td}_L(D_N)\to\infty, \quad \|N\|\to\infty.$$
\label{lem:DirichletUCTD}
\end{lemma}
{\bf Proof.} Let $N>L,$ $N=(N_1,\dots,N_d).$ Since
$\|D_N\|^2_{\td} = \prod_{j=1}^d (2 N_j + 1)$ and
$\langle{\cal B}_L D_N,D_N\rangle=0$,
$\langle{\cal A}_L D_N,D_N\rangle =~\prod_{j=1}^d (2 N_j + 1 - L_j)$ and 
\begin{equation*}
    \begin{split}
        \|{\cal B}_L D_N\|^2_{\td} 
        & = \sum\limits_{-N \le k \le N} \left( \sum_{j=1}^d (L_j k_j)^2 + \sum_{j=1}^d \sum_{n=1,n\neq j}^d L_j L_n k_j k_n  \right) \\
        & = \prod_{j=1}^d (2 N_j + 1)   \sum_{j=1}^d L_j^2 \frac{N_j(N_j+1)}{3},
    \end{split}
\end{equation*}
we obtain
{\small
\begin{equation*}
    \begin{split}
\mathrm{UP}^{\td}_L(D_N) & = \frac{1}{\|L\|^4} \left( \frac{\prod_{j=1}^d (2 N_j + 1)^2}{\prod_{j=1}^d (2 N_j + 1 - L_j)^2} -1\right) \sum_{j=1}^d L_j^2 \frac{N_j(N_j+1)}{3} \\
& = \frac{1}{\|L\|^4} \frac{\left( 1- \prod_{j=1}^d \left(1 - \frac{L_j}{2 N_j + 1}\right) \right)  \left( 1 + 
\prod_{j=1}^d \left(1 - \frac{L_j}{2 N_j + 1}\right) \right)}{\prod_{j=1}^d \left(1 - \frac{L_j}{2 N_j + 1}\right)^2} \sum_{j=1}^d L_j^2 \frac{N_j(N_j+1)}{3} 
    \end{split}
\end{equation*}
\begin{equation*}
    \begin{split}
\phantom{\mathrm{UP}^{\td}_L(D_N)} & \ge \frac{1}{\|L\|^4} \left( 1-  \left(1 - \min_{j} \frac{L_j}{2 N_j + 1}\right)^d \right)   \sum_{j=1}^d L_j^2 \frac{N_j(N_j+1)}{3} 
\\
& \ge \frac{d}{2^{d-1}\|L\|^4}  \min_{j} \frac{L_j}{2 N_j + 1}  \sum_{j=1}^d L_j^2 \frac{N_j(N_j+1)}{3},
    \end{split}
\end{equation*}}
where the last inequality is due to the mean value theorem.
Thus, $\mathrm{UP}^{\td}_L(D_N)\to\infty$ as $\|N\|\to\infty$. $\Diamond$

For the case of the multivariate Fej\'er kernel
$$
F_n(x) = \sum_{k\in\zd, \|k\|_\infty<n} \left(1-\frac{\|k\|_\infty}{n} \right)\mathrm{e}^{2\pi \mathrm{i} \langle k, x\rangle}
$$
the computation of the directional uncertainty product is  more involved.
For computations we need  the following notations
$${\cal F}(d,n)= \sum_{j=1}^n j^d, \quad  
{\cal F}^o(d,n) = \sum_{j=0}^{n-1} (2j+1)^d = {\cal F}(d,2n) -  2^d {\cal F}(d,n).
$$
Also, we need the rate of growth of the above functions when $n\to\infty.$ Due to the Faulhaber formula, we get
$$
{\cal F}(d,n)= \frac{n^{d+1}}{d+1} + \frac{n^d}{2} + \frac{d n^{d-1}}{12} + O(n^{d-3}),$$
$${\cal F}^o(d,n) = \frac{2^d n^{d+1}}{d+1} - \frac{2^{d-1} d\ n^{d-1}}{12} +  O(n^{d-3}),
$$
\begin{equation*}
{\cal F}^o(d,n-1)  = \frac{2^d}{d+1} n^{d+1} - 2^d n^d + \frac{2^d \ 11 d}{12} n^{d-1}  + O(n^{d-2}).
\end{equation*}

\begin{lemma}
    Let $F_n$ be the Fej\'er kernel, $n\in{\mathbb N},$ $L\in\zd.$ Then
    $$
    \mathrm{UP}_L^{\td}(F_n) \to \frac{(d+1)^2 (d+2)^2}{6 d \left(d+3\right)(d+4)}, \quad n\to\infty.
    $$
    \label{lem:FejerUCTD}
\end{lemma}
{\bf Proof.}
Firstly, we compute 
$$
\|F_n\|^2_{\td} = \sum_{\|k\|_\infty<n} \left( 1 - \frac{\|k\|_\infty}{n}\right)^2 = 1 + \sum_{j=1}^n \sum_{\|k\|_\infty=j} \left( 1 - \frac{j}{n}\right)^2.
$$
It is not hard to see that the number of vectors $k\in\zd$ such that $\|k\|_\infty=j$ is equal to $(2j+1)^d-(2j-1)^d$. Applying the above equalities we can estimate the rate of growth
\begin{equation*}
    \begin{split}
    \|F_n\|^2_{\td} & = 1 + \sum_{j=1}^n ((2j+1)^d-(2j-1)^d) \left( 1 - \frac{j}{n}\right)^2  
= \frac{2{\cal F}^o(d,n)}{n} -\frac{ {\cal F}^o(d+1,n)}{n^2}\\
& = \frac{2^{d+1}}{(d+1)(d+2)} n^d + \frac{2^d}{12} n^{d-2} +  O(n^{d-4}).
   \end{split}
\end{equation*}
Now we compute $\|{\cal B}_L F_n\|_{\td}^2.$
Firstly, consider $\frac{\partial F_n}{\partial x_1}.$
Let $D_n(u)$ be the one-dimensio\-nal Dirichlet kernel. Since $F_n(x) = \frac{1}{n} \sum_{j=0}^{n-1} \prod_{l=1}^d D_{j}(x_l)$, we obtain
$$
\frac{\partial F_n}{\partial x_1} = \frac{1}{n} \sum_{j=0}^{n-1} \prod_{l=2}^d D_{j}(x_l) (D_{j}(x_1))'_{x_1}.
$$
Therefore, 
$$
\left\|{\cal B}_{e_1} F_n\right\|_{\td}^2 = \frac{1}{n^2} \int_{\td} \left| \sum_{j=0}^{n-1} \prod_{l=2}^d D_{j}(x_l) \frac{(D_{j}(x_1))'_{x_1}}{2\pi} \right|^2 \mathrm{d}x 
$$
\begin{multline*}
    = \frac{1}{n^2}  \sum_{j=0}^{n-1}  \prod_{l=2}^d \int_{{\mathbb T}}  |D_{j}(x_l)|^2 \mathrm{d}x_l   \int_{{\mathbb T}} \left|\frac{1}{2\pi}\frac{\mathrm{d} D_{j}(x_1)}{\mathrm{d} x_1} \right|^2 \mathrm{d}x_1  \\ + \frac{1}{n^2}  \sum_{j=0}^{n-1} \sum_{m=0,\atop m\neq j}^{n-1} \prod_{l=2}^d \int\limits_{{\mathbb T}}  D_{j}(x_l)\overline{D_{m}(x_l)} \mathrm{d}x_l   \int\limits_{{\mathbb T}} \frac{1}{4\pi^2}\frac{\mathrm{d} D_{j}(x_1)}{\mathrm{d} x_1}\overline{\frac{ \mathrm{d} D_{m}(x_1)}{\mathrm{d} x_1}}  \mathrm{d}x_1 =: S_1 + S_2.
\end{multline*}
The first sum is equal to
$$
 S_1 = \frac{1}{n^2}  \sum_{j=0}^{n-1} (2j+1)^{d-1} \frac{j (j+1) (2j+1)}{3} =  \frac{1}{n^2}  \sum_{j=0}^{n-1} (2j+1)^{d} \frac{j (j+1)}{3}.
$$
The second sum we rewrite as follows
\begin{equation*}
    \begin{split}
S_2 & = \frac{2}{n^2}  \sum_{j=0}^{n-1} \sum_{m=j+1}^{n-1} \prod_{l=2}^d \int_{{\mathbb T}}  D_{j}(x_l)\overline{D_{m}(x_l)} \mathrm{d}x_l   \int_{{\mathbb T}} \frac{1}{4\pi^2}\frac{\mathrm{d} D_{j}(x_1)}{\mathrm{d} x_1}\overline{\frac{\mathrm{d} D_{m}(x_1)}{\mathrm{d} x_1}}  \mathrm{d}x_1 \\
& =\frac{2}{n^2}  \sum_{j=0}^{n-1} \sum_{m=j+1}^{n-1} (2j+1)^{d-1} \frac{j (j+1) (2j+1)}{3}\\ &=
\frac{2}{n^2}  \sum_{j=0}^{n-2} (n-1-j) (2j+1)^{d} \frac{j (j+1)}{3}.
   \end{split}
\end{equation*}
Combining the two sums together we get
$$
\left\|{\cal B}_{e_1} F_n\right\|_{\td}^2  = S_1 + S_2 = 
\frac{1}{n^2}  \sum_{j=0}^{n-2} (2n-1-2j) (2j+1)^{d} \frac{j (j+1)}{3} + \frac{(n-1) (2n-1)^d}{3n} .
$$
Let $c_k=c_k(F_n)$ be the Fourier coefficients of $F_n$.
Then 
$$
\|{\cal B}_L F_n\|_{\td}^2 = \sum_{k\in\zd} \langle L, k\rangle^2 c_k^2 =   \sum_{k\in\zd} \left(\sum_{j=1}^d L_j k_j \right)^2 c_k^2  
$$
$$
=\sum_{j=1}^d L^2_j \sum_{k\in\zd}  k^2_j c_k^2  + \sum_{j=1}^d \sum_{m=1,\atop m\neq j}^d  L_j L_m \sum_{k\in\zd}  k_j k_m c_k^2 
= \|L\|^2 \left\|{\cal B}_{e_1} F_n\right\|_{\td}^2  
$$
due to the symmetry of the  coefficients.
Now we establish the rate of growth  of 
{\small
$$
\|{\cal B}_L F_n\|_{\td}^2 = \frac{\|L\|^2}{3n} \left( 2 \sum_{j=0}^{n-2} (2j+1)^{d} (j^2+j) -  \sum_{j=0}^{n-2} \frac{(2j+1)^{d+1} (j^2+j)}{n}  +  (n-1) (2n-1)^d \right). 
$$}
Denote 
$G(d,n-1) = \sum_{j=0}^{n-2} (2j+1)^{d} (j^2+j)$.
It can be stated that
$$
G(d,n-1) = \frac{1}{4}\left({\cal F}^o(d+2,n-1)-{\cal F}^o(d,n-1)\right).
$$
Again, we need to estimate the rate of growth of $G(d,n-1)$ and $G(d+1,n-1).$
Thus,
$$
G(d,n-1) = \frac{1}{4}\left(\frac{2^{d+2}}{d+3} (n-1)^{d+3}  + O((n-1)^{d+2}) \right)=
\frac{2^{d}}{d+3} n^{d+3} + O( n^{d+2}).
$$
Therefore,
$$
G(d+1,n-1) = \frac{2^{d+1}}{d+4} n^{d+4} + O(n^{d+3}).
$$
Also note that $(n-1) (2n-1)^d = O(n^{d+1}).$
So, finally,
\begin{equation*}
    \begin{split}
\|{\cal B}_L F_n\|_{\td}^2 & = \frac{\|L\|^2}{3n} \left( 2 G(d,n-1)  - \frac{1}{n} G(d+1,n-1)  +  (n-1) (2n-1)^d \right)  \\
  & =\frac{\|L\|^2}{3n} \left( \frac{2^{d+1}}{d+3} n^{d+3}  - \frac{2^{d+1}}{d+4} n^{d+3}+ O(n^{d+2}) \right) \\ & = 
\frac{\|L\|^2}{3} \left( \frac{2^{d+1}}{(d+3)(d+4)} n^{d+2} + O(n^{d+1}) \right).
   \end{split}
\end{equation*}
Thus, the rate of growth of $\|{\cal B}_L F_n\|_{\td}^2 / \|F_n\|_{\td}^2$ is given by
$$
\frac{\|{\cal B}_L F_n\|_{\td}^2}{\|F_n\|_{\td}^2} = \frac{\|L\|^2}{3}  \frac{ \frac{2^{d+1}}{(d+3)(d+4)} n^{d+2} + O(n^{d+1}) }{\frac{2^{d+1}}{(d+1)(d+2)} n^d + O(n^{d-2})} = \frac{\|L\|^2}{3} \frac{(d+1)(d+2)}{(d+3)(d+4)} n^2 + O(n).
$$
Since $F_n$ is even and $\frac{\partial F_n}{\partial x_l}$, $l=1,\dots,d$, are odd, then $\langle{\cal B}_L F_n, F_n\rangle = 0$ and also $\text{var}_L^F(f) = \|{\cal B}_L F_n\|_{\td}^2 / \|F_n\|_{\td}^2.$

Now, we compute the commutator 
\begin{equation*}
  \begin{split}
\langle{\cal A}_L F_n, F_n\rangle &= \int_{\td} \mathrm{e}^{2\pi \mathrm{i} \langle L, x\rangle} F_n^2(x) \mathrm{d}x = \frac{1}{n^2} \int_{\td}  \mathrm{e}^{2\pi \mathrm{i} \langle L, x\rangle} \left| \sum_{j=0}^{n-1} \prod_{l=1}^d D_j(x_l)\right|^2 \mathrm{d}x 
\\  & = \frac{1}{n^2}  \sum_{j=0}^{n-1}  \prod_{l=1}^d \int\limits_{\mathbb T}   \mathrm{e}^{2\pi \mathrm{i} L_l x_l} \left|  D_j(x_l)\right|^2 \mathrm{d}x_l  \\ & \phantom{   }+
    \frac{1}{n^2}  \sum_{j=0}^{n-1}  \sum_{m=0, m\neq j}^{n-1}  \prod_{l=1}^d \int\limits_{\mathbb T}   \mathrm{e}^{2\pi \mathrm{i} L_l x_l} D_j(x_l)\overline{D_m(x_l)} \mathrm{d}x_l  := R_1 + R_2.
 \end{split}
\end{equation*}
Let us consider the first sum. The inner integral is the dot product of two Dirichlet kernels which are the same but one of them is shifted by $L_l$. Thus, this integral is equal to $2j+1-|L_l|.$ Hence, for big enough $n$ we get
$$
R_1 = \frac{1}{n^2}  \sum_{j=n^*}^{n-1}  \prod_{l=1}^d (2j+1-|L_l|),
$$
where $n^*$ is such that $2n^*+1-|L_l|>0$ for all $l=1,\dots,d$ and $2n^*-1-|L_l|<0$ for some $l=1,\dots,d.$ In fact, we need to compute the rate of growth of $R_1$. Applying Vieta's formulas for $R_1$ and the formulas for ${\cal F}^o(d,n)$ and ${\cal F}^o(d-1,n)$, we obtain
\begin{equation*}
  \begin{split}
R_1 & =  \frac{1}{n^2}  \sum_{j=n^*}^{n-1} \left( (2j+1)^d - \|L\|_1 (2j+1)^{d-1} \right) + O(n^{d-3}) \\ & = \frac{2^d}{d+1} n^{d-1}- \|L\|_1 \frac{2^{d-1}}{d} n^{d-2} + O(n^{d-3}),
 \end{split}
\end{equation*}
as $n\to\infty.$ 
Now, we consider  
$$
R_2 = \frac{2}{n^2}  \sum_{j=0}^{n-2}  \sum_{m=j}^{n-1}  \prod_{l=1}^d \int_{\mathbb T}   \mathrm{e}^{2\pi \mathrm{i} L_l x_l} D_j(x_l)\overline{D_m(x_l)} \mathrm{d}x_l.
$$
The inner integral is the dot product of two Dirichlet kernels which are of different size and one of them is shifted by $L_l$. Its value is equal to
$$
\int_{\mathbb T}   \mathrm{e}^{2\pi \mathrm{i} L_l x_l} D_j(x_l)\overline{D_m(x_l)} \mathrm{d}x_l  
=  2j+1 - (|L_l| - (m-j))_+.$$
Changing the variable of summation $m$ to $\w m = m-j$ in $R_2$, we obtain
$$
R_2 =  \frac{2}{n^2}  \sum_{j=0}^{n-2}  \sum_{\w m=1}^{n-1-j}  \prod_{l=1}^d ( 2j+1 - (|L_l| - \w m)_+).
$$
Applying Vieta's formulas, we get
{\small
\begin{equation*}
\begin{split}
    R_2 
& = \frac{2}{n^2}  \sum_{j=0}^{n-2}  (n-1-j) (2j+1)^d -  \frac{2}{n^2}  \sum_{j=0}^{n-2} (2j+1)^{d-1} \sum_{\w m=1}^{n-1-j} \sum_{l=1}^d (|L_l| - \w m)_+ + O(n^{d-3}),
\end{split}
\end{equation*}}
as $n\to \infty$. Note that
$$
\sum_{\w m=1}^{n-1-j} \sum_{l=1}^d (|L_l| - \w m)_+ =  \sum_{l=1}^d \sum_{\w m=1}^{|L_l|} (|L_l| - \w m) =  \sum_{l=1}^d  \frac{|L_l|(|L_l| - 1)}{2} = \frac{\|L\|^2 - \|L\|_1}{2}.
$$
Therefore,
{\small $$
R_2 = \frac{2n-1}{n^2} \sum_{j=0}^{n-2} (2j+1)^d - \frac{1}{n^2} \sum_{j=0}^{n-2} (2j+1)^{d+1}  -  \frac{2}{n^2}  \sum_{j=0}^{n-2} (2j+1)^{d-1} \frac{\|L\|^2 - \|L\|_1}{2} + O(n^{d-3}), 
$$}
as $n\to \infty$. Applying the Faulhaber formulas for ${\cal F}^o (d,n-1)$, ${\cal F}^o (d+1,n-1)$, ${\cal F}^o (d-1,n-1)$, we get
{\small\begin{equation*}
\begin{split}
R_2 & = \frac{2n-1}{n^2} {\cal F}^o (d,n-1) - \frac{1}{n^2} {\cal F}^o (d+1,n-1)  -  \frac{2}{n^2}  {\cal F}^o (d-1,n-1)\frac{\|L\|^2 - \|L\|_1}{2} \\ & =
\frac{2^{d+1}}{(d+1)(d+2)} n^d - \frac{2^d}{d+1} n^{d-1} + \left( \frac{2^d}{12} - 2^{d-1}\frac{\|L\|^2 - \|L\|_1}{d} \right) n^{d-2} + O(n^{d-3}),
\end{split}
\end{equation*}}
as $n\to \infty$. Thus, 
$$
\langle{\cal A}_L F_n, F_n\rangle = R_1 + R_2 =  \frac{2^{d+1}}{(d+1)(d+2)} n^d  + \left( \frac{2^d}{12} - \frac{2^{d-1} \|L\|^2}{d} \right)  n^{d-2} + O(n^{d-3}),
$$
as $n\to \infty$. Combining these estimates we obtain for the angular variance 
$$
\left( \frac{\|F_n\|_{\td}^2}{\langle{\cal A}_L F_n, F_n\rangle}\right)^2 -1  = \frac{\|L\|^2}{n^2}\frac{(d+1) (d+2)}{2d} + O\left(\frac{1}{n^3}\right), \quad n\to \infty.
$$
So,  the directional uncertainty product of the sequence of Fej\'er kernels is given by
\begin{equation*}
\begin{split}
\mathrm{UP}^{\td}_L(F_n) 
& = \frac{(d+1)^2 (d+2)^2}{6 d \left(d+3\right)(d+4)} + O\left(\frac{1}{n}\right) , \quad n\to\infty. \Diamond
\end{split}
\end{equation*}
For $d=1$, the limit is equal to $\frac{3}{10}$ which  coincides with the known results.
For $d=2$, the limit is equal to $\frac{2}{5}$. For $d=3$, the limit is $\frac{100}{189}.$
Similar results are valid for $\mathrm{UP}^{\td}_{GG}$.

\begin{lemma}
 Let $
D_N(x)= \sum\limits_{-N \le k \le N} \mathrm{e}^{2\pi \mathrm{i} \langle k, x\rangle},$
where $N\in \zd,$ $N>\nul,$ $L\in\zd.$ 
    Let $F_n$ be the Fej\'er kernel, $n\in{\mathbb N},$ $L\in\zd.$ Then
    $$ \mathrm{UP}^{\td}_{GG}(D_N)\to\infty, \quad \|N\|\to\infty, \quad \text{ and }\quad
    \mathrm{UP}_{GG}^{\td}(F_n) \to \frac{(d+1)^2 (d+2)^2}{6 d \left(d+3\right)(d+4)}, 
    $$
    as $n\to\infty.$
    \label{lem:DFUCTD}
\end{lemma}

{\bf Proof.}
Since $\|D_N\|^2_{\td} = \prod_{j=1}^d (2N_j+1),$  $\langle {\cal A}_{j} D_N, D_N \rangle = 2 N_j \prod_{i=1, i\neq j}^d (2N_i+1),$
$$
\|{\cal B}_{j} D_N\|^2_{\td} =\frac{N_j(N_j+1)(2 N_j+1)}{3} \prod_{i=1, i\neq j}^d (2N_i+1) = \|D_N\|^2_{\td} \frac{N_j(N_j+1)}{3},
$$
$$
\langle {\cal B}_{j} D_N, D_N \rangle = 0, \quad \frac{\langle {\cal A}_{j} D_N, D_N \rangle}{\|D_N\|^2_{\td}} = 1-\frac{1}{2N_j+1},
$$
$$
\mathrm{UP}^{\td}_{GG}(D_N) = \frac{d - \sum_{j=1}^d \left( 1-\frac{1}{2N_j+1}\right)^2}{\left(d - \sum_{j=1}^d \frac{1}{2N_j+1}\right)^2} \sum_{j=1}^d \frac{N_j(N_j+1)}{3}.
$$
Thus, $\mathrm{UP}^{\td}_{GG}(D_N)\to\infty$ as $\|N\|\to\infty.$

Concerning the Fej\'er kernel, using the rates of growths and decays established in the previous lemma, we get for $j=1,\dots,d$:  $\langle{\cal B}_j F_n, F_n\rangle = 0$,
$$
\frac{\|{\cal B}_{j} F_n\|^2_{\td}}{\|F_n\|^2_{\td}} = \frac{1}{3} \frac{(d+1)(d+2)}{(d+3)(d+4)} n^2 + O(n), \quad n\to\infty,
$$
$$
\frac{|\langle {\cal A}_{j} F_n, F_n \rangle|}{\|F_n\|^2_{\td}} = 1-\frac{(d+1) (d+2)}{4 d n^2} + O(1/n^3), \quad n\to\infty,
$$
$$
\text{var}^A_{GG}(F_n) 
= \frac{(d+1) (d+2)}{2 d^2 n^2} + O(1/n^3), \quad n\to\infty,
$$
then
$$
\mathrm{UP}^{\td}_{GG}(F_n) = \text{var}^A_{GG}(F_n)  \sum_{j=1}^d \frac{\|{\cal B}_{j} F_n\|^2_{\td}}{\|F_n\|^2_{\td}} \to \frac{(d+1)^2 (d+2)^2}{6 d \left(d+3\right)(d+4)}, \quad n\to\infty. \Diamond
$$

Also, we can place the Dirichlet and Fej\'er kernels along the direction vector $L$. Namely, let
$$
D_n^L(x) = \sum_{m=-n}^n \mathrm{e}^{2\pi \mathrm{i} \langle k_0+Lm, x\rangle},
\quad
F_n^L(x) = \sum_{m=-n}^n \left(1-\frac{|m|}{n}\right) \mathrm{e}^{2\pi \mathrm{i} \langle k_0+Lm, x\rangle},
$$ for some $k_0\in\zd.$

\begin{lemma}
    Let $L\in\zd.$ Then
    $
    \mathrm{UP}_L^{\td}(D_n^L) \to \infty, \quad \mathrm{UP}_L^{\td}(F_n^L) \to \frac{3}{10},$ as $n\to\infty.
    $
\end{lemma}

The proof can be done by straightforward computations.

\subsection{The minimal angular variance}   
   
Now we give a multivariate analogue of Rauhut's result in~\cite{Rauhut2005UC} on minimizing the angular variance for trigonometric polynomials. For a finite subset $S$ in $\zd$, denote the set of trigonometric polynomials $$\Pi_S = \left\{\sum_{k\in S} c_k \mathrm{e}^{2 \pi \mathrm{i} \langle k, x\rangle}: c_k\in {\mathbb C}\right\}.$$
Then, one is interested in best localized polynomials, i.e., for a fixed $L\in\zd$ find all trigonometric polynomials $p$, whose coefficient support is inside some fixed set $S\subset\zd$ and its directional uncertainty product takes its minimal value, i.e.
$$
    \min\limits_{ p\in \Pi_S}\{\mathrm{UP}^{\td}_L(p)\}.
$$
This problem is difficult for an arbitrary set $S.$ 

Nevertheless,  it is possible to minimize the angular frequency and the frequency variance separately. For the frequency variance the minimum value is equal to zero and it attains on trigonometric polynomials that have only one non-zero coefficient as it was shown above.

For the angular variance the situation is not so trivial. Again, since $\text{var}_A^L(f) = \text{var}_A^L(a f)$ for any appropriate $f\in L_2(\td)$ and $a\in\rd,$ we can assume that $\|p\|_{\td} = 1$.  
So, let us consider the problem
\begin{equation}
        \min\limits_{p\in \Pi_S}\{\text{var}_L^A(p): \|p\|_{\td}=1\}.
        \label{eq:MinAng}
\end{equation}
Since the set $\{p\in \Pi_S, \|p\|_{\td}=1\}$ is a compact set and $\text{var}_A^L(p)$ is continuous (except the cases when $\langle {\cal A}_L p,p \rangle=0$), we can conclude that the minimum exists.
 During the proof of the theorem below, we need to split the set $S$ into several disjoint "threads" of points. Each "thread" $U$ is a subset of $S$ that looks as follows (the order of elements is fixed)
$$
U = \{ k, k+L, k+2L, \dots, k + m L\},
$$ 
where $m \in \n $ and $k\in S$ are  chosen such that $k-L \notin S,$ and $k + (m+1) L\notin S.$ 
These "threads" are sorted by decreasing number of elements. Assume that the number of "threads" is $u$ and $U_0$ is the longest (if there are several of them, we can take any). Therefore, 
$$
S = \bigcup\limits_{i=0}^{u-1} U_i = \bigcup_{i=0}^{u-1} \{k_i, k_i+L, \dots, k_i + m_i L \}.
$$
The next Theorem states that the minimal angular variance~(\ref{eq:MinAng}) depends on the length of the longest "thread" inside $S$. 

\begin{theorem}
The minimal angular variance for trigonometric polynomials with coefficient support inside $S$ is equal to
$$\min\limits_{p\in \Pi_S}\{\text{\emph{var}}^A_L(p), \|p\|_{\td}=1\} = \tan^2 \frac{\pi}{m_0+2},$$
where $m_0+1$ is the length of the longest "thread" inside $S.$
The minimum is attained by the trigonometric polynomial $p^{\mathrm{min}}$ whose non-zero coefficients are placed on this "thread". The Fourier coefficients $c_k=c_k(p^{\mathrm{min}})$ of such a polynomial are defined as follows
\begin{equation}
    c_k = \begin{cases} \sin \frac{\pi j}{m_0+2} &\mbox{if } k = k_0 + (j-1) L, \quad j=1,\dots,m_0+1, \\
    0 & \mbox{else.} \end{cases} 
 \label{eq:minVarACoef}
\end{equation}
The directional uncertainty product is given by
$$
\mathrm{UP}^{\td}_L(p^{\mathrm{min}})  = \frac{m_0(m_0+4)}{12} \tan^2 \frac{\pi}{m_0+2} - \frac{1}{2}.
$$
\label{theo:MinAngL}
\end{theorem}

{\bf Proof.}
Note that $\text{var}_L^A(p) = \left|\sum_{k\in\zd} c_{k-L}\overline{c_k}\right|^{-2} - 1.
$
Then the minimization problem~(\ref{eq:MinAng}) is equivalent to 
$
    \max\left\{\left|\sum\limits_{k\in S} c_{k-L}\overline{c_k}\right|^2: \sum\limits_{k\in S} |c_k|^2 =1\right\}.
$
Firstly, we reduce the problem to real coefficients. Let $c_k = r_k \mathrm{e}^{i \phi_k},$ $k\in S.$ Thus, we have to maximize
$$
    \left|\sum\limits_{k\in S} r_{k-L}\overline{r_k} \mathrm{e}^{i (\phi_{k-L}-\phi_{k})} \right|^2, \quad \text{ as } \sum\limits_{k\in S} |r_k|^2 =1.
$$
The maximum is attained only if $ \mathrm{e}^{i (\phi_{k-L}-\phi_{k})} = \mathrm{const}, \forall k\in S$ or
$(\phi_{k-L}-\phi_{k}) \equiv \alpha\mod 2\pi$, $\forall k\in S$, for some $\alpha\in{\mathbb R}.$
Then we can take phases as follows $\phi_k = \beta + \alpha \frac{\langle L,k \rangle}{\|L\|^2},$ where $\beta\in{\mathbb R}.$

Therefore, the minimization problem~(\ref{eq:MinAng}) is reduced to the following
\begin{equation}
 \frac{\sum\limits_{k\in S} c_{k-L}c_k}{\sum\limits_{k\in S} c_k^2} \to \max, \quad c_k \in \rn, \quad c_k \ge 0.
 \label{eq:minQF}
\end{equation}

Let us rewrite the problem using quadratic forms. We enumerate all coefficients using one index according to the order of "threads" in $S$ and the order of the elements inside "threads". 
Hence,~(\ref{eq:minQF}) can be written in matrix form
$$
\frac{\sum\limits_{k\in S} c_{k-L}c_k}{\sum\limits_{k\in S} c_k^2} = \frac{C^\mathrm{T} M C}{C^\mathrm{T} C}, 
$$
where  $C = \{ c_k\}_{k\in S}$ is a column vector and $M$ is a block diagonal matrix
$$
M = \left(
\begin{array}{cccc}
 M_0 & 0  & \dots & 0 \\
 0 & M_1 & \dots & 0 \\
 \vdots &  \vdots & \ddots & \vdots \\
 0 & 0 & \dots &  M_{u-1} \\
\end{array}
\right), \quad
M_i=\left(
\begin{array}{ccccc}
 0 & \frac{1}{2} & 0 & \dots  & 0 \\
 \frac{1}{2} & 0 & \frac{1}{2} & \dots  & 0 \\
 0 & \frac{1}{2} & 0 & \dots  & 0 \\
  \vdots &  \vdots &  \vdots & \ddots &  \vdots \\
 0 & 0 & \frac{1}{2} & 0 & \frac{1}{2} \\
 0 & 0 & 0 & \frac{1}{2} & 0 \\
\end{array}
\right), \quad i=0,\dots,u-1.
$$
Here $M_i$ is a $(m_i+1)\times(m_i+1)$ tridiagonal Toeplitz matrix with zeros on the main diagonal and halves on the sub- and super-diagonal.
Therefore, it remains to find the maximal eigenvalue of the matrix $M$ and the corresponding eigenvector, since
$\frac{C^\mathrm{T} M C}{C^\mathrm{T} C} \le \lambda_{\mathrm{max}}(M).$
The eigenvalues of these matrices $M_i$ are known (see, e.g.,~\cite[p. 53]{Meinardus}).
They are equal to $\cos \frac{\pi n}{m_i+2}$, $n=1,\dots,m_i+1$. Since the set of eigenvalues of the block-diagonal matrix $M$ is the union of eigenvalues of its blocks, the maximum eigenvalue of $M$ is equal to $\lambda_{\mathrm{max}}(M) = \cos \frac{\pi}{m_0+2}$.
Moreover, the corresponding eigenvector also can be found. 
For the matrix $M_i$ the eigenvector $v^{(i,n)}$ corresponding to the eigenvalue $\cos \frac{\pi n}{m_i+2}$ is given coordinate-wise as follows
$$
v^{(i,n)}_j= \sin \frac{\pi n j}{m_i+2}, \quad j=1,\dots,m_i+1, \quad n=1,\dots,m_i+1.
$$
Therefore, the eigenvectors of the block-diagonal matrix $M$ can be easily defined. Hence, the eigenvector $C_{\mathrm{max}}$ corresponding to the maximal eigenvalue is given  by~(\ref{eq:minVarACoef}).
The above considerations yield that
$$\min\limits_{p\in \Pi_S}\{\text{var}^A_L(p), \|p\|_{\td}=1\} = 
\frac{1}{ \lambda^2_{\mathrm{max}}(M)} - 1 = \frac{1}{ \cos^2 \frac{\pi}{m_0+2}} - 1 = \tan^2 \frac{\pi}{m_0+2}.$$

Now we compute the directional uncertainty product for the polynomials with the minimal angular variance. In fact, it remains to compute the frequency variance:
$$
\|p^{\mathrm{min}}\|^2_{\td} = \sum\limits_{k\in S} c_k^2 = \sum\limits_{n=1}^{m_0+1} \sin^2 \frac{\pi n}{m_0+2} = \frac{m_0+2}{2},
$$
{\small
$$
\sum\limits_{k\in S} \langle L, k\rangle c_k^2 =  \sum\limits_{n=1}^{m_0+1} \langle L, k_0+nL\rangle \sin^2 \frac{\pi n}{m_0+2} = \langle L, k_0\rangle \frac{m_0+2}{2} + \|L\|^2 \sum\limits_{n=1}^{m_0+1}  n \sin^2 \frac{\pi n}{m_0+2},
$$
\begin{multline*}
        \sum\limits_{k\in S} \langle L, k\rangle^2 c_k^2  =  \sum\limits_{n=1}^{m_0+1} \langle L, k_0+nL\rangle^2  \sin^2 \frac{\pi n}{m_0+2}  \\ = \langle L, k_0\rangle^2 \frac{m_0+2}{2} + 2 \langle L, k_0\rangle \|L\|^2 \sum\limits_{n=1}^{m_0+1}  n \sin^2 \frac{\pi n}{m_0+2} + \|L\|^4 \sum\limits_{n=1}^{m_0+1}  n^2 \sin^2 \frac{\pi n}{m_0+2}.
\end{multline*}}
Based on trigonometric formulas, the formulas for the Dirichlet kernel and for the conjugate Dirichlet kernel, and taking the derivatives of those kernels we compute
\begin{equation*}
    \begin{split}
    \text{var}^F_L(p^{\mathrm{min}}) & =  \frac{\sum\limits_{k\in\zd} 
\langle L, k\rangle^2 |c_k|^2}{\|p^{\mathrm{min}}\|^2_{\td}} - \left(\frac{\sum\limits_{k\in\zd} \langle L, k\rangle |c_k|^2}{\|p^{\mathrm{min}}\|^2_{\td}} \right)^2 \\ &= \|L\|^4  \left( \frac{m_0^2+4 m_0}{12} - \frac{1}{2} \cot^2 \frac{\pi}{m_0+2}\right).
    \end{split}
\end{equation*}
Finally,
$$
\mathrm{UP}^{\td}_L(p^{\mathrm{min}})  = \frac{m_0(m_0+4)}{12} \tan^2 \frac{\pi}{m_0+2} - \frac{1}{2}. 
$$
This finishes the proof. $\Diamond
$

Note, for $m_0\to\infty$, we obtain $\mathrm{UP}^{\td}_L(p^{\mathrm{min}}) \to \frac{\pi^2}{12}-\frac{1}{2} \approx 0.3224 > \frac 14.$ 

Next, we establish a similar result for the uncertainty product defined by Goh and Goodman in case where the coefficients support $S$ is a rectangle $S = \prod_{j=1}^d [-N_j,N_j] \cap \zd,$ where all $N_j>0$. The problem is to minimize
$$
\text{var}_{GG}^A(p) =  \left(d -\sum\limits_{j=1}^d \left|\sum\limits_{k\in S} c_{k-e_j}\overline{c_k}\right|^2\right) \left(\sum\limits_{j=1}^d\left|\sum\limits_{k\in S} c_{k-e_j}\overline{c_k}\right|\right)^{-2} \quad 
$$
when $\|p\|^2_{{\mathbb T}^d}=1.$ From the above considerations it follows that for $\|p\|^2_{{\mathbb T}^d}=1$ the sum
$\left|\sum\limits_{k\in S} c_{k-e_j}{c_k}\right|$ for any $j=1,\dots,d$ cannot be greater than $\cos \frac{\pi}{m+1},$ where $m$ is the length of the longest "thread". Since $S = \prod_{j=1}^d [-N_j,N_j] \cap \zd$ and $L=e_j$, then $\left|\sum\limits_{k\in S} c_{k-e_j}{c_k}\right|\le \cos \frac{\pi}{2N_j+2}$, $ j=1,\dots,d$ and for fixed $j$ the equality $\left|\sum\limits_{k\in S} c_{k-e_j}{c_k}\right| =  \cos \frac{\pi}{2N_j+2}$ attains, if
{\small
$$
 c_k =    \begin{cases} \frac{1}{\sqrt{N_j+1}}\sin \frac{\pi l}{2N_j+2} &\mbox{if } k_j = -N_j-1 + l, \ \  l=1,\dots,2N_j+1, \quad k_i=0, \text{ for } i\neq j,\\
    0 & \mbox{else,} \end{cases} 
$$}
where $k=(k_1,\dots,k_d).$

If it is possible to achieve for some $p_{\mathrm {min}}$ those values $ \cos \frac{\pi}{2N_j+2}$ for all $j=1,\dots,d$ simultaneously, then we get the minimal possible value for $\text{var}_{GG}^A(p)$ which is equal to 
\begin{equation}
\text{var}_{GG}^A(p_{\mathrm {min}}) = \left(d  - \sum_{j=1}^d \cos^2 \frac{\pi}{2N_j+2}\right) \left( \sum_{j=1}^d \cos \frac{\pi}{2N_j+2} \right)^{-2}.
\label{eq:VarAGGMin}
\end{equation}


\begin{theorem}
The value $\min\limits_{p\in \Pi_S}\{\text{\emph{var}}^A_{GG}(p), \|p\|_{\td}=1\}$ is given by~(\ref{eq:VarAGGMin}) and it is attained by $p_{\mathrm{min}}$ if the Fourier coefficients $c_k$ of $p_{\mathrm{min}}$  are given by
    $$
    c_k = \prod_{j=1}^d \frac{1}{\sqrt{N_j+1}} \sin \frac{\pi l_j}{2 N_j + 2},
    $$
    where $l_j=k_j+N_j+1$, $j=1,\dots,d$, $k\in S.$
    \label{theo:MinAngGG}
\end{theorem}

{\bf Proof.} Let us show that the maximum values  $\left|\sum\limits_{k\in S} c_{k-e_j}{c_k}\right| =  \cos \frac{\pi}{2N_j+2}$ are attained for all $j=1,\dots,d$ simultaneously. This can be checked by direct computations. Let us fix $i=1,\dots,d$. Therefore,
$$
\sum\limits_{k\in S} c_{k-e_i}{c_k} = \prod_{j=1}^d \frac{1}{N_j+1} \sum\limits_{l_1=1}^{2N_1+1} \dots\sum\limits_{l_d=1}^{2N_d+1} \prod_{j=1}^d \sin \frac{\pi l_j}{2 N_j + 2} \sin \frac{\pi (l_j-\delta_{ij})}{2 N_j + 2}
$$
$$
=\prod_{j=1}^d \frac{1}{N_j+1} \left(
\prod_{j=1, j\neq i}^d \sum\limits_{l_j=1}^{2N_j+1} \sin^2 \frac{\pi l_j}{2 N_j + 2} \right) \sum\limits_{l_i=1}^{2N_i+1}  \sin \frac{\pi l_i}{2 N_i + 2} \sin \frac{\pi (l_i-1)}{2 N_i + 2}.
$$
It remains to note that $\sum\limits_{l_j=1}^{2N_j+1} \sin^2 \frac{\pi l_j}{2 N_j + 2} = N_j+1$ and 
{\small $$
 \sum\limits_{l_i=1}^{2N_i+1}  \sin \frac{\pi l_i}{2 N_i + 2} \sin \frac{\pi (l_i-1)}{2 N_i + 2} = 
\frac{1}{2} \sum\limits_{l_i=1}^{2N_i+1}  \sin \frac{\pi l_i}{2 N_i + 2} \left(  \sin \frac{\pi (l_i-1)}{2 N_i + 2} +  \sin \frac{\pi (l_i+1)}{2 N_i + 2}\right)  
$$
$$
=\sum\limits_{l_i=1}^{2N_i+1}  \sin^2 \frac{\pi l_i}{2 N_i + 2} \cos \frac{\pi}{2N_i+2} = (N_i+1) \cos \frac{\pi}{2N_i+2}. \Diamond
$$}

§§§§§§§§§§§§§§§§§§§§§§§§§§§§§§§§§§§§§§§§

\section{Well localized multivariate periodic Parseval wavelet frames}

First of all, we recall the notion of a Parseval frame. 
Let $H$ be a separable Hilbert space. If there exist  constants $A,\,B>0$
such that for any $f \in H$ the following inequality holds
$$
A \|f\|^2 \leq \sum_{n=1}^{\infty} \left|(f,\,f_n)\right|^2 \leq B \|f\|^2,
$$ 
then the sequence $(f_n)_{n \in \mathbb{N}}$ is called  a frame  for $H.$
A frame is a complete system. Moreover, any element $f\in H$ can be expanded in a series 
$\sum_n\alpha_n f_n,$ 
$\alpha_n\in \cn,
$ 
with respect to a frame. However, the series expansion is not unique.  
If $A=B(=1),$ then the sequence $(f_n)_{n \in \mathbb{N}}$ is called a tight frame (a Parseval frame) for $H.$

In this section we design a family of well-localized multivariate periodic Parseval wavelet frames. This is a generalization of the wavelet family constructed in \cite{LebPres14}. It turns out that these wavelet frames have optimal localization with respect to the dimension $d$ of the torus $\td$. More precisely, we claim that
$$
\lim_{j \to \infty} \mathrm{UP}_L^{\mathbb{T}^d}(\psi_{j}) =
\frac{1}{4}\frac{(d+2)(d^2-2d+4)}{d^3},
$$
so 
$
\lim_{d \to \infty}\lim_{j \to \infty} \mathrm{UP}_L^{\mathbb{T}^d}(\psi_{j})=\frac14.
$

Let
$A \in \z^{d\times d}$ be a dilation matrix that means that all the eigenvalues of the matrix are greater then $1$. The determinant of $A$ is equal to $2$. Therefore,    a full collection of coset representatives of $\mathbb{Z}^d/A \mathbb{Z}^d$ consists of $2$ elements  (see, e.g.,~\cite{KPS}). We denote these collection as $\{0,\,k_0\}$. Further,  
$
B = A^\mathrm{T},
$
$
K_j=\mathbb{Z}^d \cap B^{j} [-1/2,\,1/2)^d.
$
Put by definition
$
f_{j}(k) = {\rm exp} \left(-\frac{\|L\|^2 \|k\|^2}{j(j-1)}\right), 
$
where
$
L \in \mathbb{Z}^d,
$ 
$j \geq 2.$ 
Let us define a $B^j$-periodic sequence  $\nu_j(k)$ 
\begin{equation}
\label{nu}
\nu_j(k) =
\left\{
\begin{array}{ll}
f_j(k) & k \in \text{int}(K_{j-1}), \\
\left(1-f_j^2(k - B^{j-1} k_0)\right)^{1/2} &  k - B^{j-1} k_0 \in \text{int}(K_{j-1}), \\
1/\sqrt{2} & k \in \overline{K}_{j-1} \setminus \text{int}(K_{j-1}), \\
\end{array}
\right.
\end{equation}
where 
$\text{int}(K_{j})=\mathbb{Z}^d \cap B^{j} (-1/2,\,1/2)^d$
and 
$\overline{K}_{j}=\mathbb{Z}^d \cap B^{j} [-1/2,\,1/2]^d.$
The $B^j$-periodicity means that
$
\nu_j(k+B^j p) = \nu_j(k)$
 for any 
$
k,p \in \mathbb{Z}^d
$.

For instance, if $d=2$, 
$
A=
 \left(
\begin{array}{cc}
1     & 1 \\
  -1   & 1
\end{array}
\right),
$
then on the main period  
$K_j$ the sequence $\nu_j(k)$ 
is defined as follows
$$
\nu_j(k) =
\left\{
\begin{array}{ll}
f_j(k) & k \in \text{int}(K_{j-1}), \\
\left(1-f_j^2(k - B^{j-1} (v_1(r) \ v_2(r))^\mathrm{T})\right)^{1/2} & B^{j-2[j/2]} k \in Q_r, k \in K_{j} \setminus K_{j-1}, \\
1/\sqrt{2} & k \in \overline{K}_{j-1} \setminus \text{int}(K_{j-1}), \\
\end{array}
\right.
$$
where
$Q_r$ is the  $r$-th quadrant of $\mathbb{R}^2$, $v_1(r)= -\cos (\pi r/2),$ $v_2(r)= -\sin (\pi r/2),$ 
$r=1,\dots,4,$ $[y]=\text{max}\{n\in \mathbb{N} \ : \ n\le y \}.$ 
Finally, let us define an auxiliary function 
$\xi_j \in L_2(\mathbb{T}^d)$ with the Fourier coefficients
$$
\widehat{\xi}_{j}(k):=\prod_{r=j+1}^{\infty}\nu_r(k). 
$$
Later, in Theorem \ref{main_loc}, we will prove that the infinite product converges. 
 Then scaling masks, scaling functions, wavelet masks, and wavelet functions are defined respectively as 
\begin{equation}
\label{uep}
\begin{array}{l}
	\mu_j(k):=\sqrt{2} \nu_j(k),\\
\widehat{\varphi_{j}}(k):=2^{-j/2}\widehat{\xi}_{j}(k), 
  \\
	\lambda_j(k):=\mathrm{e}^{2\pi \mathrm{i} \langle k_0, B^{-j}k\rangle }\mu_j(k+B^{j-1} k_0), \\
	{\widehat{\psi}}_{j}(k):=\lambda_{j+1}(k) \widehat{\varphi}_{j+1}(k).
\end{array}
\end{equation}

\begin{theorem}
\label{main_loc}
Suppose $\varphi_j,$  $\psi_j$   are the functions defined in (\ref{uep}) and $\nu_j$ is a sequence defined in (\ref{nu}). Then the set 
$
\Psi = \{\varphi_j,\, \psi_j(\cdot - (A^{-j} k)\}_{j \in \mathbb{N}\cup\{0\},k \in \mathcal{L}_j},
$
where
$\mathcal{L}_j$ is a full collection of coset representatives of $\mathbb{Z}^d/A^j \mathbb{Z}^d,$ 
forms a Parseval frame of $L_2(\mathbb{T}^d)$, 
and the following equalities hold true 
\begin{equation}
\label{lim_2}
\lim_{j \to \infty} \mathrm{UP}_L^{\mathbb{T}^d}(\varphi_{j}) = 1/4, 
\quad
\lim_{j \to \infty} \mathrm{UP}_L^{\mathbb{T}^d}(\psi_{j}) = \frac{1}{4}\frac{(d+2)(d^2-2d+4)}{d^3}. 
\end{equation}
 \end{theorem}

The scheme of the proof repeats  in the main features  Theorem 4 \cite{LebPres14}.  At the same time, there are differences concerning  technical details. In particular, we have to provide a new proof for an analogue of Lemma 3 \cite{LebPres14} since the existing proof can not be rewritten for the multivariate case.
We exploit Lemma 2  \cite{LebPres14}, so we cite it here for convenience.

\begin{lemma}[Lemma 2 \cite{LebPres14}]
\label{tech}
Suppose $\alpha, \beta, \gamma \in \mathbb{R},$ $m=0,1,\dots,$ and $\ 0<\,b\,<b_1,$
where $b_1$ is an absolute constant, then 
\begin{multline*}
\sum_{k \in\mathbb{Z}}(\alpha k^2 + \beta k + \gamma)^m\, \mathrm{e}^{-b (\alpha k^2 +
 \beta k + \gamma)} \\ = (-1)^m
\frac{\partial^m}{\partial b^m}\left(
\mathrm{exp}\left(-b\left(\gamma-\frac{\beta^2}{4\alpha}\right)\right)
\sqrt{\frac{\pi}{b \alpha}}
\right)
+\mathrm{exp}\left(-\frac{\pi^2-\varepsilon}{b\alpha}\right)\,O(1),
\end{multline*}
as $b\to 0,$
where $\varepsilon>0$ is an arbitrary small parameter.
\end{lemma}

We need also several technical lemmas.

\begin{lemma}
\label{exp_ball}
Let $B^{\theta}_j(0) \subset \rd$ be a ball centered at the origin with  radius $\frac{1}{2}(1+\theta)^j,$ $\theta > 0.$ Let $M \in \z^{d\times d}$ be a dilation matrix with determinant equal to $2$. Then there exists $\theta_0>0$ and $j_0 \in \mathbb{N}$ such that 
$B^{\theta_0}_j(0) \subset M^j\td$ for $j\geq j_0.$
\end{lemma}

{\bf Proof.}
Let $\rho = \inf_j\|M^{-j}\|^{1/j}$ be the spectral radius of the matrix $M$. Since $M$ is a dilation matrix, it follows that  $\rho<1.$ Given $0<\varepsilon<1-\rho$, there exists $j_0\in \mathbb{N}$ such that $(\rho+\varepsilon)^j\geq \|M^{-j}\|$ for $j \geq j_0.$ 
Suppose $x\in B^{\theta}_j(0)$, that is $\|x\|\leq \frac{1}{2}(1+\theta)^j$,   
then  
$$
\|M^{-j}x\|\leq \|M^{-j}\| \|x\|  \leq \frac{1}{2}(1+\theta)^j (\rho+\varepsilon)^j \mbox{ for }
j\geq j_0.
$$
Therefore, any $\theta$ satisfying the inequality 
$0<\theta<(\rho+\varepsilon)^{-1}-1$ can be chosen as $\theta_0$.
This concludes the proof of Lemma \ref{exp_ball}. $\Diamond$

\begin{lemma}
\label{dif_int_sum}
Suppose 
$b=b(h)=2 h^2/(1-h),$ $0<h<1,$  and
{\small $$
F(x):=
b \|L\|^2 \|x\| \|x-L\| \left(
1-\frac{1}{4} b \|L\|^2 \left(\|x\|^2+\|x-L\|^2
\right)\right)
{\rm exp} (-h\|L\|^2 (\|x\|^2+\|x-L\|^2)),
$$}
then
$$
\sum_{k \in \zd}F(k) = \int_{ \rd}F(x)\,\mathrm{d}x + O(h^{2}) 
\mbox{ as }
h \to 0.
$$
\end{lemma}

{\bf Proof.}
 The Poisson summation formula
$$
\sum_{k \in \zd}F(k) = \sum_{k \in \zd}\widehat{F}(k) 
$$
shows that it is sufficient to prove 
$$
\sum_{k \in \zd \setminus\{0\}}\widehat{F}(k) = O(h^{2})
\mbox{ as }
h \to 0.
$$
So, we need only to find the Fourier transform of $F$.
To this end, we rewrite the function $F$ as
$$
F(x)=\frac{2}{1-h} f_1(x)f_1(x-L) - 
\frac{1}{(1-h)^2} f_2(x)f_1(x-L) -
\frac{1}{(1-h)^2} f_1(x)f_2(x-L),
$$
where
$
f_1(x)=h \|L\| \|x\| {\rm exp}(-h\|L\|^2 \|x\|^2),
$
$
f_2(x)=h^3 \|L\|^3 \|x\|^3 {\rm exp}(-h\|L\|^2 \|x\|^2).
$
Therefore, $\widehat{F}$ can be written as 
$$
\widehat{F}(\xi)=\frac{2}{1-h} \widehat{f_1}\ast \widehat{f_1(\cdot -L)}(\xi) - 
\frac{1}{(1-h)^2} \widehat{f_2}\ast \widehat{f_1(\cdot -L)}(\xi) -
\frac{1}{(1-h)^2} \widehat{f_1}\ast \widehat{f_2(\cdot -L)}(\xi).
$$
It follows from elementary properties of the Fourier transform that 
$$
\widehat{f_1}(\xi) = \|L\|^{-d} h^{1/2-d/2}\widehat{f_3}\left(\frac{\xi}{h^{1/2}\|L\|}\right), 
\quad
\widehat{f_2}(\xi) = \|L\|^{-d} h^{3/2-d/2}\widehat{f_4}\left(\frac{\xi}{h^{1/2}\|L\|}\right), 
$$
where
$$
f_3(x) = \|x\| {\rm exp}(- \|x\|^2),
\quad
f_4(x) = \|x\|^3 {\rm exp}(- \|x\|^2).
$$
Since $f_3$ is a radial function, we can exploit Theorem 3.3 chapter IV \cite{SW}. So, we get
$$
\widehat{f_3}(\xi) = 2 \pi \|\xi\|^{-\frac{d}{2}+1} \int_0^{\infty} r^{\frac{d}{2}+1} 
\mathrm{e}^{-r^2} J_{\frac{d-2}{2}}(2 \pi \|\xi\| r) \, \mathrm{d}r,
$$
where $J_n$ is a Bessel function of the first kind.
By \cite[Formula 11.4.28, p. 486]{AS} we conclude
$$
\widehat{f_3}(\xi) =
\pi^{d/2} \frac{\Gamma(d/2+1/2)}{\Gamma(d/2)} M(d/2+1/2,\, d/2,\, -\pi^2 \|\xi\|^2),
$$
where
$M$ is Kummer's (confluent hypergeometric) function. The asymptotic  behavior as $\xi \to \infty$ of this function is known and can be found, for instance, in  \cite[Formula 13.1.05, p. 504]{AS}, therefore we obtain
$$
\widehat{f_3}(\xi) = -1/2\, \pi^{-d/2-3/2} \Gamma(d/2+1/2) \|\xi\|^{-1-d}(1+O(\|\xi\|^{-2}))
\mbox{ as } \xi \to \infty.
$$
Analogously 
$$
\widehat{f_4}(\xi) = 3/4\, \pi^{-d/2-7/2} \Gamma(d/2+3/2) \|\xi\|^{-3-d}(1+O(\|\xi\|^{-2}))
\mbox{ as } \xi \to \infty.
$$
Thus, we get
$$
\widehat{f_1}(\xi) = C_1(d) h \|\xi\|^{-1-d} (1+O(h \|\xi\|^{-2}))
\mbox{ as } \xi \to \infty,
 $$
 $$
\widehat{f_2}(\xi) = C_2(d)  h^3 \|\xi\|^{-3-d} (1+O(h \|\xi\|^{-2}))
\mbox{ as } \xi \to \infty,
$$
and
$
C_1(d)=-1/2 \, \|L\| \pi^{-d/2-3/2} \Gamma(d/2+1/2),
$
$
C_2(d)=3/4 \, \|L\|^3 
\pi^{-d/2-7/2} \Gamma(d/2+3/2).
$
Since, in addition, the functions $\widehat{f_1}$ and $\widehat{f_2}$
are bounded, the convolutions 
$
\widehat{f_1}\ast \widehat{f_1(\cdot -L)}(\xi),
$
$
\widehat{f_2}\ast \widehat{f_1(\cdot -L)}(\xi), 
$
and
$
\widehat{f_1}\ast \widehat{f_2(\cdot -L)}(\xi)
$
are well-defined. Therefore,
$$
\widehat{F}(k)=O(\widehat{f_1}\ast \widehat{f_1(\cdot -L)}(k)) = O(h^2\|k\|^{-2}) 
\mbox{ as } h \to 0.
$$
Thus,
$$
\sum_{k \in \zd \setminus\{0\}}\widehat{F}(k) = O(h^{2})
\mbox{ as }
h \to 0,
$$
which  proves the result. $\Diamond$

\begin{lemma}
\label{xi0}
Let $\xi^0_j \in L_2(\td)$ be a function defined by its Fourier coefficients
$$
\widehat{\xi^0_j}(k) = {\rm exp} (-\|L\|^2 \|k\|^2/j).
$$
Then
$$
\lim_{j \to \infty} \mathrm{UP}_L^{\mathbb{T}^d}(\xi^0_{j}) = \frac{1}{4}.
$$
\end{lemma}

{\bf Proof.}
Since the coefficient $\widehat{\xi^0_j}(k)$
is a product of $d$ one-dimensional corresponding coefficients, it follows that it is sufficient to apply Lemma \ref{tech} to the series $\sum_{k \in \zd}\left|\widehat{\xi^0_j}(k)\right|^2,$ 
$\sum_{k \in \zd} \langle L, k\rangle^2 \left|\widehat{\xi^0_j}(k)\right|^2,$ and $\sum_{k \in \zd}\widehat{\xi^0_j}(k-L) \overline{\widehat{\xi^0_j}(k)},$ and then to substitute the results  
into the definition of $\mathrm{UP}_L^{\td}$. The result follows by simple computations. $\Diamond$

\begin{lemma}
\label{eta}
Let $\eta_j \in L_2(\td)$ be a function defined by its Fourier coefficients
$$
\widehat{\eta_j}(k) = \mathrm{e}^{2\pi \mathrm{i} \langle k_0, B^{-j}k\rangle }\left(1-{\rm exp}\left(-\frac{2\|L\|^2 \|k\|^2}{j(j+1)}\right)\right)^{1/2} {\rm exp} \left(-\frac{\|L\|^2 \|k\|^2}{j+1}
\right).
$$
Then
$$
\lim_{j \to \infty} \mathrm{UP}_L^{\mathbb{T}^d}(\eta_{j}) = \frac{1}{4}\frac{(d+2)(d^2-2d+4)}{d^3}.
$$
\end{lemma}

{\bf Proof.}
Denote 
$h=1/(j+1),$  $b=b(h) = 2h^2/(1-h) = 2/(j(j+1)).$
To estimate 
 $\sum\limits_{k\in\zd} |\widehat{\eta_j}(k)|^2,$
$\sum\limits_{k\in\zd} \langle L, k\rangle^2 |\widehat{\eta_j}(k)|^2$
one can use Lemma \ref{tech}
as it was described in Lemma 3 \cite{LebPres14}. 
Namely, 
$$
\sum\limits_{k\in\mathbb{Z}^d} |\widehat{\eta_j}(k)|^2 =
\sum\limits_{k\in\mathbb{Z}^d} \left( {\rm exp}(- 2 h \|L\|^2 \|k\|^2) - {\rm exp}\left(- \frac{2 h}{1-h} \|L\|^2 \|k\|^2\right)  \right) 
$$
$$
=
\prod_{n=1}^d\sum\limits_{k_n\in\mathbb{Z}}{\rm exp}(- 2 h \|L\|^2 k_n^2) 
-\prod_{n=1}^d\sum\limits_{k_n\in\mathbb{Z}}
{\rm exp}\left(- \frac{2 h}{1-h} \|L\|^2 k_n^2\right)
$$
$$
= 
\left(\frac{\pi}{2 h \|L\|^2}\right)^{d/2} - \left(\frac{\pi (1-h)}{2 h \|L\|^2}\right)^{d/2}
+O(\mathrm{e}^{-h^{-1}}).
$$
Since 
$\sum_{k\in \mathbb{Z}^d}k_n k_m |\widehat{\eta_j}(k)|^2=0,$
we analogously get
\begin{equation*}
    \begin{split}
    \sum\limits_{k\in \mathbb{Z}^d} \langle L, k\rangle^2 |\widehat{\eta_j}(k)|^2
& =
\|L\|^2 
\sum\limits_{k\in \mathbb{Z}^d}  k_1^2 |\widehat{\eta_j}(k)|^2 \\
& = 
\|L\|^2 
\frac{ \pi^{d/2}}{2(2 h \|L\|^2)^{d/2+1}}\left(1-(1-h)^{d/2+1}\right)
+O(\mathrm{e}^{-h^{-1}}),
    \end{split}
\end{equation*}
where 
$k=(k_1, k_2,\dots,k_d)^\mathrm{T}.$

However, 
we have to provide an alternative way to estimate  
$
\sum\limits_{k\in\zd} \widehat{\eta_j}(k-L) \overline{\widehat{\eta_j}(k)}.
$
We write
$$
\left|\sum\limits_{k\in\mathbb{Z}^d} \widehat{\eta_j}(k-L) \overline{\widehat{\eta_j}(k)}\right|
=
\sum_{k\in\mathbb{Z}^d} \left(1-{\rm exp}\left(-b\|L\|^2 \|k\|^2\right)\right)^{1/2}
$$
$$
\times
\left(1-{\rm exp}\left(-b\|L\|^2 \|k-L\|^2\right)\right)^{1/2}
{\rm exp} (-h\|L\|^2 (\|k\|^2+\|k-L\|^2)).
$$
Using the Taylor formula for the function 
$f(b)=(1-\mathrm{exp}(-b \|L\|^2\|k\|^2))^{1/2}(1-\mathrm{exp}(-b \|L\|^2\|k-L\|^2))^{1/2}$ 
in the neighborhood of $b=0$, we get
$$
f(b)=b \|L\|^2 \|k\| \|k-L\| \left(
1-\, b \|L\|^2 \left(\|k\|^2+\|k-L\|^2
\right)/4\right)+ f'''(\bar{d})\, b^3/6,
$$
and
$f'''(\bar{d})b^3=O(\|k\|^6 h^6).$
The last equality is deduced in Lemma 3 \cite{LebPres14}. 
So, using Lemma \ref{tech} for the remainder of the series we get 
{\small $$
\sum_{k\in\mathbb{Z}^d}
\frac{f'''(\bar{d})}{6} b^3
{\rm exp} (-h\|L\|^2 (\|k\|^2+\|k-L\|^2))
=
O\left(h^6\sum_{k\in\mathbb{Z}^d}\|k\|^6{\rm exp} (-h\|k\|^2)\right)
=O(h^{3-d/2}).
$$}
Therefore, 
$$
\left|\sum\limits_{k\in\mathbb{Z}^d} \widehat{\eta_j}(k-L) \overline{\widehat{\eta_j}(k)}\right|
=
b \|L\|^2 \sum\limits_{k\in\mathbb{Z}^d} \|k\| \, \|k-L\| \left(
1-\frac{1}{4} b \|L\|^2 \left(\|k\|^2+\|k-L\|^2
\right)\right)
$$
$$
\times
{\rm exp} (-h\|L\|^2 (\|k\|^2+\|k-L\|^2))
+O(h^{3-d/2}).
$$

Next, by Lemma \ref{dif_int_sum} we replace the series by the integral   
{\small $$
b \|L\|^2 \int\limits_{\mathbb{R}^d} \|x\| \, \|x-L\| \left(
1-\frac{1}{4} b \|L\|^2 \left(\|x\|^2+\|x-L\|^2
\right)\right)
{\rm exp} (-h\|L\|^2 (\|x\|^2+\|x+L\|^2))
\, \mathrm{d}x,
$$}
change the variable $y=x-L/2$ to obtain
{\small $$
b \|L\|^2 \int\limits_{\mathbb{R}^d}
\left\|y-\frac{L}{2}\right\| \left\|y+\frac{L}{2}\right\| 
\left(
1-\frac{1}{4} b \|L\|^2 \left(2 \|y\|^2+\frac{\|L\|^2}{2}
\right)\right)
{\rm exp} \left(-h\|L\|^2 \left( 2 \|y\|^2+\frac{\|L\|^2}{2}\right)\right)\,\mathrm{d}y,
$$}
and convert it to the polar coordinates $y=r \beta,$ where
$$
\beta = (\sin \phi_1 \sin \phi_2 \dots \sin \phi_{d-1},\,\cos \phi_1 \sin \phi_2 \dots \sin \phi_{d-1},\dots ,\cos \phi_{d-1})^\mathrm{T}.
$$
So, we get
{\small 
$$
b \|L\|^2 {\rm exp}\left(-h\frac{\|L\|^4}{2}\right) \int\limits_{[0,\,\infty)\times[0,\,2\pi)
\times[0,\,\pi)^{d-2}} r^2
\left(
\left(1+r^{-2}\frac{\|L\|^2}{4}\right)^2 - r^{-2} \langle L, \beta\rangle^2
\right)^{1/2}
$$
$$
\times
\left(
1-\frac{1}{4} b \|L\|^2 \left(2 r^2+\frac{\|L\|^2}{2}
\right)\right)
{\rm exp} (-2 h \|L\|^2  r^2)
$$
$$
\times r^{d-1}
\sin \phi_2 \sin^2 \phi_3 \dots \sin^{d-2} \phi_{d-1}
\,\mathrm{d}r\,\mathrm{d}\phi.
$$}
Next,   applying the Taylor formula 
for 
$\left(
(1+r^{-2}\|L\|^2/4)^2 - r^{-2} \langle L, \beta\rangle^2
\right)^{1/2}$ 
with respect to $1/r$, changing the variable $(2 h)^{1/2} \|L\|  r = t$,
integrating with respect to $\phi$,
and recalling that $b=2h^2(1-h)^{-1}$,
we obtain 
{\small \begin{multline*}
\frac{2h}{1-h} {\rm exp}\left(-h\frac{\|L\|^4}{2}\right)
\left(\frac{\pi}{2h\|L\|^2}\right)^{d/2} 
\frac{1}{\Gamma(d/2)} \\ \times
\int_0^{\infty}t^{d+1} {\rm exp} (-t^2) 
\left(1+\frac{d-2}{2d} \frac{\|L\|^4}{t^2} h\right)
\left(1-\frac{h t^2}{2(1-h)}\right)\,\mathrm{d}t.
\end{multline*}}
Integrating with respect to $t$, we finally obtain
{\small $$
\left|\sum\limits_{k\in\mathbb{Z}^d} \widehat{\eta_j}(k-L) \overline{\widehat{\eta_j}(k)}\right|
$$
$$
=
\frac{h}{1-h} {\rm exp}\left(-h\frac{\|L\|^4}{2}\right) \left(\frac{\pi}{2h\|L\|^2}\right)^{d/2} 
\left(\frac{d}{2}+
\frac{d-2}{2d} \|L\|^4 h-
\frac{h}{1-h} \frac{d(d+2)}{8}  +O(h^2)\right).
$$}
It is easy to see that 
$\sum\limits_{k\in\zd} \langle L, k\rangle |\widehat{\eta_j}(k)|^2=0.$
It remains to substitute the expressions for 
$\sum\limits_{k\in\zd} |\widehat{\eta_j}(k)|^2,$
$\sum\limits_{k\in\zd} \langle L, k\rangle^2 |\widehat{\eta_j}(k)|^2,$
and
$\sum\limits_{k\in\zd} \widehat{\eta_j}(k-L) \overline{\widehat{\eta_j}(k)}$
to the definition of $\mathrm{UP}_L^{\td}$. Lemma \ref{eta} is proved. $\Diamond$

\vspace{0.4cm}

{\bf Proof of Theorem \ref{main_loc}.}
First of all, it is straightforward to see that the infinite product
$
\widehat{\xi}_{j}(k):=\prod_{r=j+1}^{\infty}\nu_r(k) 
$
converges. As usual, if an infinite product is equal to zero then it is also considered convergent. Indeed, it follows from (\ref{nu}) that 
$\nu_j(k)=f_j(k)$ for $k \in {\rm int} (K_{j-1})$, and Lemma \ref{exp_ball} 
says that there exist $\theta_0>0$ and $j_0 \in \mathbb{N}$ such that $\|k\|\leq (1+\theta_0)^j/2$ 
implies $k \in {\rm int} (K_{j-1})$ for $j\geq j_0$. Therefore, 
$\nu_j(k)=f_j(k)$ for $\|k\|\leq (1+\theta_0)^j/2$. So, we get
\begin{equation}
\label{xi}
\widehat{\xi}_{j}(k) =
\left\{
\begin{array}{ll}
\prod\limits_{r=j+1}^{j_1}\nu_r(k)
\prod\limits_{r=j_1+1}^{\infty} f_r(k) =
\left(\prod\limits_{r=j+1}^{j_1}\nu_r(k)\right)
{\rm exp}\left(-\frac{\|L\|^2 \|k\|^2}{j_1}\right), &
j < j_1 \\
\prod\limits_{r=j+1}^{\infty} f_r(k)= {\rm exp}\left(-\frac{\|L\|^2 \|k\|^2}{j} \right), &
j\geq j_1
\end{array}
\right.
\end{equation}
where $j_1 = \lfloor\log_{1+\theta_0}(2\|k\|)\rfloor+1.$  Therefore, 
$\widehat{\xi}_{j}(k)$ is well-defined and $\xi_j \in L_2(\td).$
Then one can check that all conditions of the unitary extension principle are fulfilled for the functions $\varphi_j,\, \psi_j$ (see Theorem 2.2 \cite{GT1}).
Therefore, the set 
$
\Psi = \{\varphi_j,\, \psi_j(\cdot - (A^{-j} k)\}_{j \in \mathbb{N}\cup\{0\},k \in \mathcal{L}_j},
$
forms a Parseval frame of $L_2(\mathbb{T}^d).$
 
To check (\ref{lim_2}),
as in the univariate case, we
 introduce two auxiliary functions 
$\xi^0_j$ and $\eta_j$ by the Fourier coefficients
$$
\widehat{\xi^0_j}(k) = {\rm exp} (-\|L\|^2 \|k\|^2/j),
$$
$$
\widehat{\eta_j}(k) = \mathrm{e}^{2\pi \mathrm{i} \langle k_0 B^{-j}k \rangle}\left(1-{\rm exp}\left(-\frac{2\|L\|^2 \|k\|^2}{j(j+1)}\right)\right)^{1/2} {\rm exp} \left(-\frac{\|L\|^2 \|k\|^2}{j+1}
\right).
$$
Now we claim that
$$
\lim_{j\to \infty}\|\xi^0_j-\xi_j\|_{L_2(\mathbb{T}^d)}+\sum_{n=1}^d\|(\xi^0_j-\xi_j)_n'\|_{L_2(\mathbb{T}^d)}=0,
$$
$$
\lim_{j\to \infty}\|\eta_j-2^{j/2}\psi_j\|_{L_2(\mathbb{T}^d)}+\sum_{n=1}^d\|(\eta_j-2^{j/2}\psi_j)_n'\|_{L_2(\mathbb{T}^d)}=0,
$$
where $f'_n$ means again the partial derivative of $f$ with respect to $x_n.$
Indeed, 
Since 
$
\widehat{\xi^0_j}(k) = \widehat{\xi_j}(k) 
$
and 
$
\widehat{\eta_j}(k) =2^{j/2} \widehat{\psi_j}(k)  
$
for $k\in {\rm int}(K_{j-1})$, and, therefore, for $\|k\|\leq (1+\theta_0)^j/2$ (see Lemma \ref{exp_ball}), it follows that 
$$
\|\xi^0_j-\xi_j\|^2_{L_2(\mathbb{T}^d)}+\sum_{n=1}^d\|(\xi^0_j-\xi_j)_n'\|^2_{L_2(\mathbb{T}^d)}
$$
$$
=
\sum_{\|k\|\geq (1+\theta_0)^j/2}\left|\widehat{\xi_j}(k) - \widehat{\xi^0_j}(k)\right|^2
+
4\pi^2\sum_{n=1}^d \sum_{\|k\|\geq (1+\theta_0)^j/2}k_n^2\left|\widehat{\xi_j}(k) - \widehat{\xi^0_j}(k)\right|^2.
$$
By (\ref{xi}), we have
$$
\left|\widehat{\xi_j}(k) - \widehat{\xi^0_j}(k)\right| \leq 
2{\rm exp}\left(-\frac{\|L\|^2 \|k\|^2}{j_1}\right).
$$
Substituting this majorant to the series, we get that the series tends to zero as $j \to \infty$ as  a remainder of a convergent series. For the functions $\eta_j$ and $2^{j/2}\psi_j$ it can be checked analogously. 
The functional $\mathrm{UP}_L^{\mathbb{T}^d}$ is continuous with respect to the norm 
$\|f\|_{L_2(\mathbb{T}^d)}+\sum_{n=1}^d\|f_n'\|_{L_2(\mathbb{T}^d)}$, which can be checked as in the one-dimensional case in Lemma 1 \cite{LebPres14}. Moreover,  
$\mathrm{UP}_L^{\mathbb{T}^d}(\xi^0_j),$ 
$\mathrm{UP}_L^{\mathbb{T}^d}(\eta_j)$ are bounded with respect to $j$, which follows from Lemma \ref{xi0} and Lemma \ref{eta}. Therefore,
$$
\lim_{j \to \infty} 
\mathrm{UP}_L^{\mathbb{T}^d}(\varphi_j) = \lim_{j \to \infty} \mathrm{UP}_L^{\mathbb{T}^d}(\xi^0_j),
\quad
\lim_{j \to \infty} 
\mathrm{UP}_L^{\mathbb{T}^d}(2^{j/2}\psi_j) =
\lim_{j \to \infty} \mathrm{UP}_L^{\mathbb{T}^d}(\eta_j).
$$
Finally, the functional $\mathrm{UP}_L^{\mathbb{T}^d}$  is  homogeneous that is 
$\mathrm{UP}_L^{\mathbb{T}^d}(\alpha f) = \mathrm{UP}_L^{\mathbb{T}^d}(f)$, $\alpha\neq 0.$ 
So
$$
\mathrm{UP}_L^{\mathbb{T}^d}(2^{j/2}\psi_j)=\mathrm{UP}_L^{\mathbb{T}^d}(\psi_j).
$$
Thus, 
$$
\lim_{j \to \infty} 
\mathrm{UP}_L^{\mathbb{T}^d}(\varphi_j) = \lim_{j \to \infty} \mathrm{UP}_L^{\mathbb{T}^d}(\xi^0_j),
\quad
\lim_{j \to \infty} 
\mathrm{UP}_L^{\mathbb{T}^d}(\psi_j) =
\lim_{j \to \infty} \mathrm{UP}_L^{\mathbb{T}^d}(\eta_j).
$$
To conclude the proof of Theorem \ref{main_loc} it remains to apply 
Lemma \ref{xi0} and Lemma \ref{eta}. $\Diamond$

\begin{acknowledgements}
The authors thank Professor O.L.Vinogradov for a valuable observation in the proof of Lemma \ref{eta}.
\end{acknowledgements}


%
%

\end{document}